\documentclass[11pt]{gtart}
\usepackage{amsmath, amssymb, graphicx, epsfig, graphs, verbatim, psfrag}
\usepackage{epic}

\newtheorem{thm}{Theorem}[section]

\newtheorem{cor}[thm]{Corollary}
\newtheorem{lem}[thm]{Lemma}
\newtheorem{prop}[thm]{Proposition}

\theoremstyle{definition}
\newtheorem{defn}[thm]{Definition}

\newtheorem{rem}[thm]{Remark}

\numberwithin{equation}{section}

\newcommand{\al}{\alpha}

\newcommand{\si}{\sigma}

\newcommand{\ze}{\zeta}

\newcommand{\J}{{\mathcal {J}}}
\newcommand{\ot}{{\overline {\t}}}
\newcommand{\os}{{\overline {\s}}}
\newcommand{\oa}{{\overline {a}}}

\newcommand{\x}{\times}
\newcommand{\s}{\mathbf s}
\renewcommand{\u}{\mathbf u}
\renewcommand{\t}{\mathbf t}

\newcommand{\Z}{\mathbb Z}
\newcommand{\N}{\mathbb N}
\newcommand{\Q}{\mathbb Q}
\newcommand{\R}{\mathbb R}

\newcommand{\CP}{{\mathbb C}{\mathbb P}}

\newcommand{\del}{\partial}

\newcommand{\lra}{\longrightarrow}
\newcommand{\hf}{{{\widehat {HF}}}}

\DeclareMathOperator{\rot}{rot}

\DeclareMathOperator{\PD}{PD}
\DeclareMathOperator{\Spin}{Spin}

\begin{document}

\title{Ozsv\'ath--Szab\'o invariants and tight\\ contact
three--manifolds, II}

\author{Paolo Lisca}
\address{Dipartimento di Matematica\\
Universit\`a di Pisa \\I-56127 Pisa, Italy} 
\email{lisca@dm.unipi.it}

\secondauthor{Andr\'{a}s I. Stipsicz}
\secondaddress{R\'enyi Institute of Mathematics\\
Hungarian Academy of Sciences\\
H-1053 Budapest\\ 
Re\'altanoda utca 13--15, Hungary \emph{and}\\
Institute for Advanced Study, Princeton, New Jersey}
\secondemail{stipsicz@math-inst.hu}

\begin{abstract}
Let $p$ and $n$ be positive integers with $p>1$, and let $E_{p,n}$ be
the oriented 3--manifold obtained by performing $p^2n-pn-1$ surgery
on a positive torus knot of type $(p, pn+1)$. We prove that $E_{2,n}$
does not carry tight contact structures for any $n$, while $E_{p,n}$
carries tight contact structures for any $n$ and any odd $p$. In
particular, we exhibit the first infinite family of closed, oriented,
irreducible 3--manifolds which do not support tight contact
structures. We obtain the nonexistence results via standard methods of
contact topology, and the existence results by using a quite delicate
computation of contact Ozsv\'ath--Szab\'o invariants.
\end{abstract}
\primaryclass{57R17} 
\secondaryclass{57R57} 
\keywords{tight contact
structures, contact surgery, Seifert fibered 3--manifolds,
Ozsv\'ath--Szab\'o invariants}

\maketitle

\section{Introduction}\label{s:intro}

Let $S^3_r(K)$, $r\in\Q$, be the oriented 3--manifold obtained by
performing rational $r$--surgery along a knot $K\subset
S^3$. In~\cite{LSuj} we used the Ozsv\'ath--Szab\'o invariants to
study the existence of tight contact structures on $S^3_r(K)$.  In
particular, we proved that if $T_{p,q}$ is the positive $(p,q)$ torus
knot, then $S^3_r(T_{p,q})$ carries positive, tight contact
structures for every $r\neq pq-p-q$.

On the other hand, it was proved by Etnyre and Honda~\cite{EH} that
$S^3_1(T_{2,3})$ supports no positive tight contact structure.
Therefore, the question whether the 3--manifolds
$S^3_{pq-p-q}(T_{p,q})$ carry positive, tight contact structures seems
to be particularly interesting.

Consider the oriented 3--manifold
\[
E_{p,n}:=S^3_{p^2n-pn-1}(T_{p,pn+1}) 
\]
The first main result of this paper is the following:

\begin{thm}\label{t:bound}
Let $p$, $n$ be positive integers with $p>1$. Then, the number of
isotopy classes of tight contact structures carried by $E_{p,n}$ is at
most 
\[
2 \max\{p(p-1)-4, 0\}.
\]
\end{thm}

An immediate corollary of Theorem~\ref{t:bound} is:

\begin{cor}\label{c:notight}
Let $n$ be a positive integer. Then, the oriented 3--manifold
$E_{2,n}$ admits no positive, tight contact structures.\qed
\end{cor}

Notice that Corollary~\ref{c:notight} generilizes the result of Etnyre
and Honda~\cite{EH}. Since the 3--manifolds $E_{2,n}$ are Seifert
fibered with base $S^2$ and three exceptional fibers, by~\cite{Wa}
they are irreducible. Therefore, Corollary~\ref{c:notight} gives the
first infinite family of closed, oriented and irreducible 3--manifolds
not carrying positive, tight contact structures.

In the second part of the paper we prove the following:

\begin{thm}\label{t:T34}
Let $n$, $p$ be positive integers with $p>1$ odd. Then, $E_{p,n}$ 
carries positive, tight contact structures.
\end{thm}

In order to motivate this result, we also prove that the oriented
3--manifolds $E_{p,n}$ do not support any fillable contact structures
(Proposition~\ref{p:nofill}). Therefore, one cannot prove the
existence of tight contact structures by presenting the 3--manifolds
$E_{p,n}$ as boundaries of symplectic fillings. In fact, we need to
use the more sophisticated methods provided by Heegaard Floer theory.

The paper is organized as follows. In Section~\ref{s:two} we prove
Theorem~\ref{t:bound} and so verify Corollary~\ref{c:notight}.  The
proof uses convex surface theory along the lines of~\cite{EH, GhS}. In
the second part of the paper (Sections~\ref{s:three} to \ref{s:six})
we prove Theorem~\ref{t:T34} using the Ozsv\'ath--Szab\'o invariants.
In Section~\ref{s:three} we recall the relevant facts of Heegaard
Floer theory. In Section~\ref{s:four} we show that the 3--manifolds
$E_{p,n}$ do not support symplectically fillable contact structures.
In Section~\ref{s:five} we define suitable contact structures on the
manifolds $E_{p,n}$ ($p>1$ odd) and in Section~\ref{s:six} we verify
their tightness. The techniques used in the first part of the paper
(Section~\ref{s:two}) are completely independent from the methods
applied in the second part
(Sections~\ref{s:three}--\ref{s:six}). However, the two approaches
nicely complement each other, in the sense that using both of them on
the same 3--manifold appears to be an effective way to attack the
classification problem for tight contact structures.

{\bf Acknowledgments.} The first author was partially supported by
  MURST, and he is a member of EDGE, Research Training Network
  HPRN-CT-2000-00101, supported by The European Human Potential
  Programme.  The second author was partially supported by OTKA T034885.
  Part of this collaboration took place when the second author visited
  the University of Pisa. He wishes to thank the Geometry Group of the
  Mathematics Department for hospitality and support.

\section{Proof of Theorem~\ref{t:bound}}
\label{s:two}

We will follow the methods developed in~\cite{EH} and implemented
in~\cite{GhS}. We will assume that the reader is familiar with the
theory of convex surfaces~\cite{Gi0} as well as the
references~\cite{EH, GhS}. 

\begin{figure}[ht]
\begin{center}
\setlength{\unitlength}{1mm}
\begin{picture}(55,20)
\put(50.5,11){$0$}
\put(11,11){$-\frac 1a$}
\put(23,10.5){$-\frac 1b$}
\put(35,11){$-\frac 1c$}
\includegraphics[height=2cm]{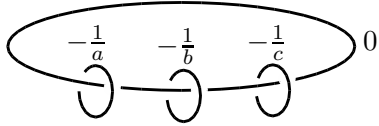}
\end{picture}    
\end{center}
\caption{\quad Surgery diagram for the Seifert 
fibered 3--manifold $M(a,b,c)$}
\label{f:seifert}
\end{figure}

We now recall the notations used in~\cite{EH, GhS}. Denote the Seifert
fibered 3--manifold given by the surgery diagram of
Figure~\ref{f:seifert} by $M(a,b,c)$ (with $a,b,c\in \Q$).

\begin{lem}\label{l:surg}
Let $p,n\in\N$ with $p\geq 2$ and $n\geq 1$. Then, there exists an
orientation--preserving diffeomorphism
\[
S^3_{p^2n-pn-1}(T_{p,pn+1})\cong
M\left(-\frac{1}{p}, \frac{n}{pn+1}, \frac{1}{p(n+1)+1}\right).
\]
\end{lem}

\begin{figure}[ht]
\begin{center}
\setlength{\unitlength}{1mm}
\begin{picture}(113,100)
\put(0,0){\includegraphics[height=10cm]{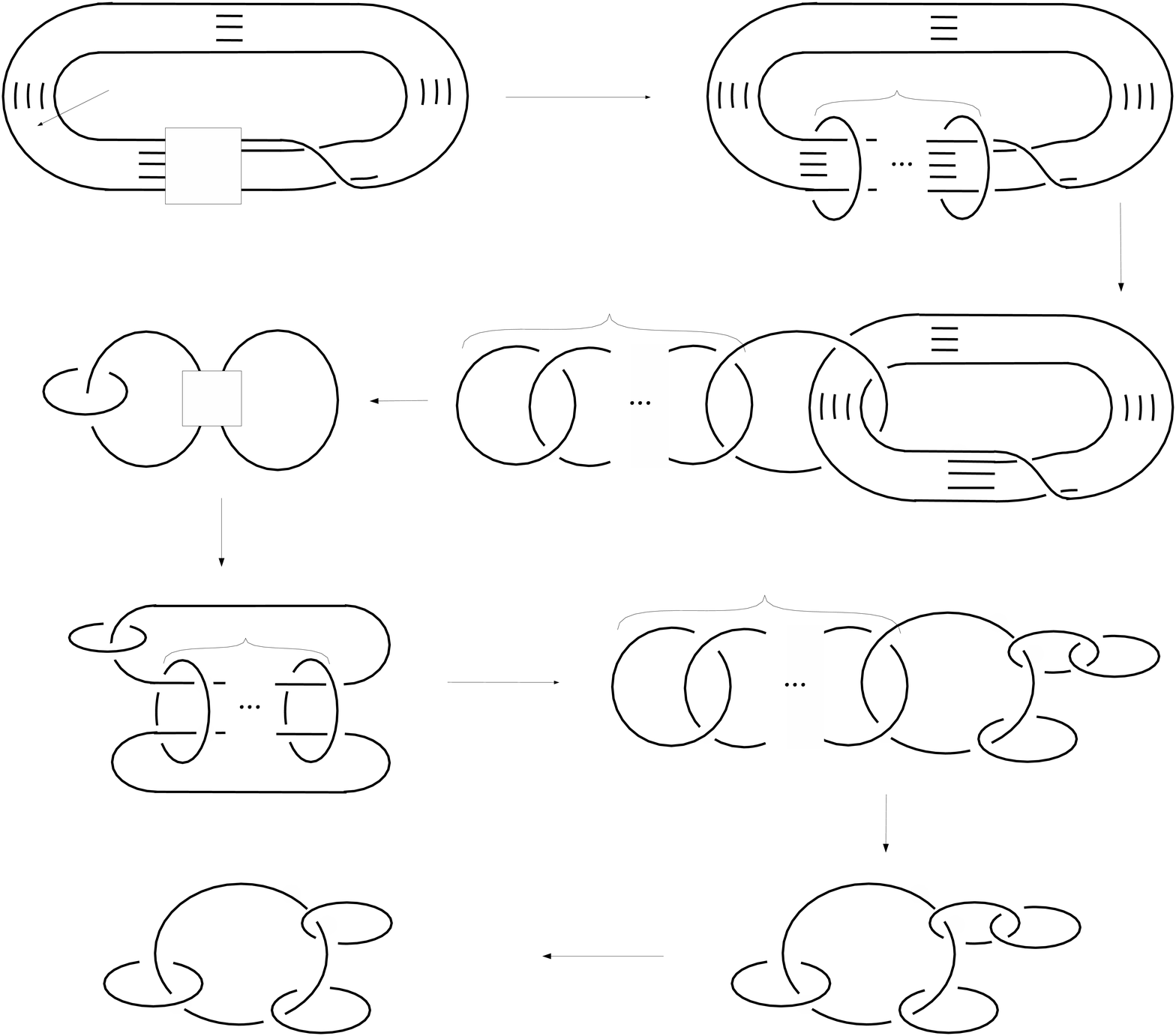}}
\put(20,99){\footnotesize $r$}
\put(12,90){\footnotesize $p$ strands}
\put(19,82){\footnotesize $n$}
\put(85,99){\footnotesize $r-np^2$}
\put(84,90){\footnotesize $n$}
\put(77,75){\footnotesize $-1$}
\put(89,75){\footnotesize $-1$}
\put(46,52){\footnotesize $-2$}
\put(53,52){\footnotesize $-2$}
\put(63,52){\footnotesize $-2$}
\put(72,51){\footnotesize $-1$}
\put(54,69.5){\footnotesize $n-1$}
\put(89,70){\footnotesize $r-np^2$}
\put(2.5,61){\footnotesize $n$}
\put(20,60){\footnotesize $p$}
\put(14,68){\footnotesize $0$}
\put(28,68){\footnotesize $r-np^2$}

\put(30,43){\footnotesize $-p$}
\put(5,38){\footnotesize $n$}
\put(23,39.5){\footnotesize $p$}
\put(11,31.5){\footnotesize $-1$}
\put(33,31.5){\footnotesize $-1$}
\put(17,22){\footnotesize $r-p(np+1)$}

\put(69,43){\footnotesize $p-1$}
\put(60,25.6){\footnotesize $-2$}
\put(67,25.6){\footnotesize $-2$}
\put(78,25.6){\footnotesize $-2$}
\put(86,42){\footnotesize $-1$}
\put(95,24){\footnotesize $r-p(np+1)$}
\put(96,40.7){\footnotesize $-p$}
\put(110,36.5){\footnotesize $n$}

\put(95,3){\footnotesize $r-p(np+1)$}
\put(102.5,11.5){\footnotesize $n$}
\put(90,15.3){\footnotesize $-p$}
\put(75,15){\footnotesize $0$}
\put(66.5,6){\footnotesize $p$}

\put(36,2){\footnotesize $r-p(np+1)$}
\put(8,6){\footnotesize $p$}
\put(33,15.3){\footnotesize $-p-\frac 1n$}
\put(16,15){\footnotesize $0$}

\end{picture}    
\end{center}
\caption{A diffeomorphism between $S^3_r(T_{p,pn+1})$ and 
$M(-\frac{1}{p}, \frac{n}{pn+1}, \frac{1}{p(np+1)-r})$}
\label{f:kirby}
\end{figure}

\begin{proof}
An orientation--preserving diffeomorphism is given by the sequence of
Kirby moves of Figure~\ref{f:kirby} for $r=p^2n-pn-1$ (see
e.g.~\cite{GS} for an introduction to Kirby calculus).
\end{proof}

Define 
\[
E_{p,n}:=S^3_{p^2n-pn-1}(T_{p,pn+1}).
\]
In view of Lemma~\ref{l:surg} and following~\cite{EH, GhS}, we start
by decomposing $E_{p,n}$ into $S^1\times \Sigma _0$, where $\Sigma _0$
is $S^2$ minus three disks, and three copies of $S^1\times D^2$
identified with neighbourhoods $V_i$ of the singular fibers $F_i$,
$i=1,2,3$.  In order to recover $E_{p,n}$ from $S^1\times \Sigma _0$
we need to glue these three copies of $S^1\times D^2$ to its three
boundary tori. We can prescribe the gluing maps by matrices once we
fix identifications of the boundary tori with $\R^2/\Z ^2$. To do
that, for each boundary component of $\partial (S^1\times \Sigma _0)$
we identify the intersection with a section $\{ * \} \times \Sigma _0$
with the image of the line $\langle (1,0)\rangle$, and the fiber with
the image of the line $\langle (0,1)\rangle$. For the boundaries of
the solid tori $S^1\times D^2$, the meridional direction is uniquely
determined by the property of being homologically trivial in
$S^1\times D^2$. The longitude is unique only up to a
$\Z$--action. This indeterminacy results in a certain degree of
freedom in choosing the particular gluing matrices. We choose:
\[
A_i\co\del(S^1\x D^2)\to -\del(E_{p,n}\setminus  V_i),\quad i=1,2,3,
\]
\[
A_1=\left(\begin{matrix}
 p & -1 \\
 1 & 0 
\end{matrix}
\right),\ \ \ 
A_2= 
\left(\begin{matrix}
 pn+1 & pn-p+1 \\
 -n & 1-n
\end{matrix}
\right),\ \ \ 
A_3=
\left(\begin{matrix}
 p(n+1)+1 & 1 \\
 -1 & 0
\end{matrix}
\right). 
\]
The matrices $A_i$ have determinant one, and the ratios of the
elements in their first columns equal the surgery coefficients
appearing in the surgery diagram. We shall denote by $F_i$ the
singular fibers inside the glued-up tori, while each neighbourhood of
$F_i$ (as a subspace of $E_{p,n}$) will be called $V_i$,
$i=1,2,3$. From the matrices $A_i$ it is immediate to compute that a
regular fiber of the fibration has slope
\[
v_1=p,\quad v_2=-\frac{pn+1}{pn-p+1}\quad\text{and}
\quad v_3=-(p(n+1)+1)
\]
when viewed, respectively, in $\del V_i$, $i=1,2,3$, while the
meridian of each $V_i$ has slope
\[
c_1=\frac{1}{p},\quad c_2=-\frac{n}{pn+1}\quad
\text{and}\quad c_3=-\frac{1}{p(n+1)+1}
\]
when viewed in $-\del(E_{p,n}\setminus  V_i)$, $i=1,2,3$. The
numbers $v_1$, $v_2$ and $v_3$ are called the~\emph{vertical slopes},
while $c_1$, $c_2$ and $c_3$ are the~\emph{critical slopes}.

Recall that the~\emph{slope} of a convex torus in standard form
identified with $\R^2/\Z^2$ is, by definition, the slope of any
component of its dividing set.

\rk{Remark.} If $T$ is a convex torus in standard form isotopic to
$\partial V_i$ and the slope of $T$ with respect to the identification
$-\del(E_{p,n}\setminus  V_i)\cong\R^2/\Z^2$ given above is equal to the
critical slope of $F_i$, then the contact structure under
consideration is overtwisted. In fact, any Legendrian divide on $T$
bounds an overtwisted disk in $V_i$.

Let $f\subset E_{p,n}$ be a Legendrian curve isotopic to a regular fiber of
the fibration. There are two framings of $f$: the one coming from the
fibration and the one induced by the contact structure. The difference
between the fibration framing and the contact framing is, by
definition, the~\emph{twisting number} of $f$.

Let $F_i$ be a Legendrian singular fiber with twisting number $m_i$
and standard neighbourhood $V_i$. Then, the slope of the torus
$\partial V_i$ is $\frac{1}{m_i}$ with respect to the identification
$\partial V_i\cong\R^2/\Z^2$ given above. The same slope is equal to,
respectively,
\[
b_1=\frac{m_1}{pm_1-1},\quad 
b_2=-\frac{n(m_2+1)-1}{(pn+1)m_2+p(n-1)+1}
\]
and
\[
b_3=-\frac{m_3}{(p(n+1)+1)m_3+1}
\]
when computed with respect to the chosen identification $-\del
(E_{p,n}\setminus V_i)\cong\R^2/\Z^2$.  The numbers $b_1$, $b_2$ and $b_3$
are called the~\emph{boundary slopes}.

\begin{lem}\label{l:twisting}
Let $\xi$ be a positive, tight contact structure on $E_{p,n}$. Then,
the singular fibers $F_1$, $F_2$ and $F_3$can be isotoped to
Legendrian positions such that
\[
m_1=0\quad\text{and}\quad m_2=m_3=-1.
\]
Moreover, we can find (nonstandard) neighbourhoods $V'_i\supset V_i$ 
with convex boundaries such that the slopes of $-\del(E_{p,n}\setminus V'_i)$
are all infinite.
\end{lem}

\begin{proof}
The argument is a simple adaptation of the proof
of~\cite[Lemma~7]{EH}. Notice that the statement
of~\cite[Lemma~7]{EH} coincides with the statement we want to prove
for $(n,p)=(1,2)$. Therefore, we will assume $(n,p)\neq (1,2)$.

Let $V_2$ and $V_3$ be standard neighbourhoods of $F_2$ and $F_3$ with
vertical rulings on their boundaries. Up to stabilizing $F_2$ and
$F_3$, we may assume $m_2,m_3<-1$. Then, there are two possible cases.

{\bf Case I.} Suppose there is a vertical annulus $A$ between $V_2$
and $V_3$ having no boundary parallel dividing curves. Then, by the
Imbalance Principle~\cite[Proposition~3.17]{H1},
\begin{equation}\label{e:m2m3}
(pn+1)m_2+p(n-1)+1=(p(n+1)+1)m_3+1,
\end{equation}
that is,
\[
m_3=\frac{(pn+1)m_2+p(n-1)}{p(n+1)+1} = 
m_2+1-\frac{pm_2+2p+1}{p(n+1)+1}.
\]
Since $m_3\in\Z$, this implies that $p(n+1)+1\geq 7$ divides
$pm_2+2p+1\neq 0$, therefore $m_2<-2$ and we have
\[
\vert pm_2+2p+1\vert =p\vert m_2 \vert -2p-1\geq p(n+1)+1.
\]
This observation implies that Equation~\eqref{e:m2m3} can hold only if
\[
\vert m_2\vert \geq n+3+ \frac{2}{p},
\]
i.e., if $m_2\leq -(n+4)$. If
we cut along $A$ and round corners, we get a torus $T$ of slope
\begin{equation}
-s_T=-\frac{n(m_2+1)+\frac{(pn+1)m_2+p(n-1)}{p(n+1)+1}}{(pn+1)m_2+p(n-1)+1} .
\end{equation}
surrounding the fibers $F_2$ and $F_3$. When viewed as minus the
boundary of the complement of a neighbourhood of $F_1$, the slope of $T$
is $s_T$. We claim that
\begin{equation}\label{e:claim}
s_T>\frac{m_1}{pm_1-1}
\end{equation}
In fact, it is easy to check that $s_T$ is a strictly decreasing function 
of $m_2$, and takes the value $s_T=\frac{1}{p}$ for 
\[
m_2=-1-\frac{p^2n}{p^2n-p-1}.
\]
Moreover, an easy calculation shows that, since $(n,p)\neq (1,2)$, 
\[
-(n+4) < -1-\frac{p^2n}{p^2n-p-1}.
\]
It follows that for $m_2\leq -(n+4)$ we have $s_T>\frac{1}{p}$. 
Therefore, since 
\[
\frac{1}{p}>\frac{1}{p-\frac{1}{m_1}}=\frac{m_1}{pm_1-1},
\]
the claim~\eqref{e:claim} is proved. This immediately implies the
existence of a convex vertical torus $T'$ with slope $\infty$. Then,
let $A_i$, $i=1,2,3$, be vertical convex annuli between a Legendrian
divide of $T'$ and a ruling of $\del V_i$, $i=1,2,3$. As long as
$m_i<0$, we can find bypasses on $A_i$ attached to $\del V_i$ for each
$i=1,2,3$. By attaching those bypasses to $V_i$ we can find bigger
standard neighbourhoods of the singular fibers $F_i$, which amounts to
increasing the twisting numbers $m_i$ as long as the assumptions of
the Twist Number Lemma~\cite[Lemma~6]{EH} hold, i.e., as long as
\[
\frac{1}{p}\geq m_1+1, \quad -\frac{pn-p+1}{pn+1}\geq m_2+1, \quad   
-\frac{1}{p(n+1)+1}\geq m_3+1.
\]
Consequently, we can increase the $m_i$'s up to $m_1=0$ and
$m_2=m_3=-1$. Moreover, the Legendrian divide of $T'$ allows us to
attach further vertical bypasses to the standard neighbourhoods until
we obtain the neighbourhoods $V'_i$ of the statement.

{\bf Case II.} Suppose there is a vertical annulus $A$ between $V_2$
and $V_3$ with some boundary parallel dividing curve. Then, we can
attach a vertical bypass to either $V_2$ or $V_3$ and increase either
$m_2$ or $m_3$. Since under Case~I we have proved the statement, we
may assume that we fall again under Case~II. Using
Equation~\eqref{e:m2m3} it is easy to check that if $m_2=-1$ we can
always attach a vertical bypass to $V_3$ as long as
$m_3<-1$, while if $m_3=-1$ we can
attach a vertical bypass to $V_2$ as long as $m_2<-1$. Therefore, we
may assume to be able to increase $m_2$ and $m_3$ until $m_2=m_3=-1$.
At this point the values of the boundary slopes $b_2$ and $b_3$ are
\[
b_2=-\frac 1p\quad\text{and}\quad
b_3=-\frac 1{p(n+1)}.
\]
We can keep attaching vertical bypasses until the slopes of the
resulting neighbourhoods are both $-\frac 1k$, for some $0\leq k\leq
p$. Since for $k=0$ this gives a vertical convex torus of infinity
slope and the conclusion follows as in Case~I, we may assume that at
some point we can find an annulus $A$ between the two neighbourhoods
with no boundary parallel curves. After cutting and rounding we get a
torus of slope $-\frac{1}{k}$ surrounding $F_2$ and $F_3$, which can
be viewed as a torus of slope $s=\frac{1}{k}$ around $V_1$. For $k=p$,
$s$ is the critical slope of the first singular fiber, hence its
existence contradicts the tightness of $\xi$. For $0\leq k<p$ we have
\[
b_1=\frac{m_1}{pm_1-1}=\frac{1}{p-\frac{1}{m_1}}<\frac{1}{k}.
\]
Therefore there is a torus of slope $\infty$ around $F_1$,
and the conclusion follows as before.
\end{proof}

\begin{figure}[ht]
\begin{center}
\setlength{\unitlength}{1mm}
\begin{picture}(120,100)
\put(0,0){\includegraphics[height=10cm]{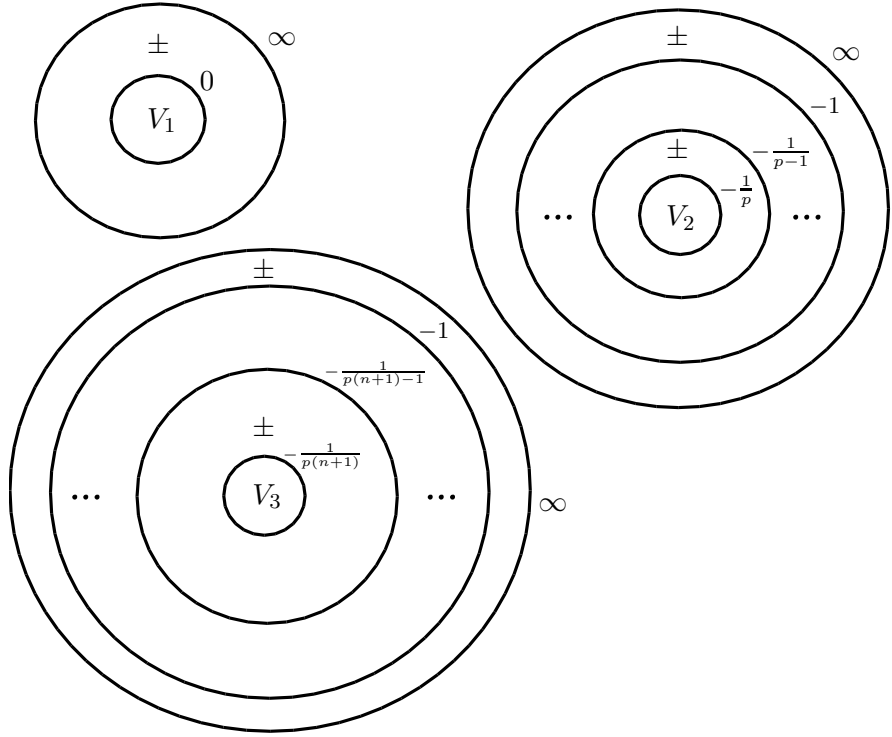}}

\put(18, 81){$V_1$}
\put(18,91){$\pm$}
\put(25,86){$0$}
\put(34,92){$\infty$}

\put(87,68){$V_2$}
\put(87,77){$\pm$}
\put(87,92){$\pm$}
\put(94,72){\footnotesize $-\frac 1p$}
\put(98.3,77){\footnotesize $-\frac 1{p-1}$}
\put(106,83){\footnotesize $-1$}
\put(109,90){$\infty$}

\put(32,31){$V_3$}
\put(32,40){$\pm$}
\put(32,61){$\pm$}
\put(36,37){\tiny $-\frac 1{p(n+1)}$}
\put(41.5,48){\tiny $-\frac 1{p(n+1)-1}$}
\put(54,53){\footnotesize $-1$}
\put(70,30){$\infty$}

\end{picture}    
\end{center}
\caption{A tight contact structure with infinity boundary slopes on
$\cup_i V'_i$}
\label{f:poss}
\end{figure}

Using Lemma~\ref{l:twisting}, we can assume the boundary slopes to be
\[
b_1=0,\quad b_2=-\frac{1}{p}\quad 
\text{and}\quad b_3=-\frac{1}{p(n+1)}.
\]
Let $V'_i$ $(i=1,2,3)$ be the neighbourhoods given in the statement of
Lemma~\ref{l:twisting}. Each of the thickened tori $V'_i\setminus V_i$ has a
decomposition into basic slices. Following the notation of~\cite{GhS}, 
any tight contact structure on $\cup_i V'_i$  with infinity boundary slopes 
can be represented and is uniquely determined by a diagram as in 
Figure~\ref{f:poss} for some choice 
of signs, where each sign denotes the corresponding type of basic slice.
Let $q_i$ denote the number of `$+$' signs in
$V_i$. Then,
\[
q_1\in \{ 0,1\},\quad 
q_2\in \{ 0, \ldots , p\}\quad\text{and}\quad 
q_3\in \{ 0, \ldots , p(n+1)\}.  
\]
Let us denote by $\xi(q_1,q_2,q_3)$ the contact structure on $\cup_i V'_i$ 
corresponding to the vector $(q_1, q_2, q_3)$.

\begin{lem}\label{l:potatos}
Let $\xi$ be a positive contact structure on $E_{p,n}$ such that  
\[
\xi|_{\cup_i V'_i} = \xi(q_1,q_2,q_3).
\] 
If $q_2\leq q_3\leq q_2+pn$, then $\xi$ is overtwisted. 
\end{lem}

\begin{proof}
By contradiction, suppose that $\xi$ is tight. The
assumption is equivalent to
\begin{equation}\label{e:potatos}
q_3\geq q_2\quad\text{and}\quad p(n+1)-q_3\geq p-q_2.
\end{equation}
Denote by $V''_2$ and $V''_3$ the neighbourhoods of $F_2$ and $F_3$,
respectively, bounded by vertical tori inside $V'_2$ and $V'_3$ with
slope $-\frac{1}{p}$. Since by~\cite[Lemma~4.14]{H1} the basic slices
of $V'_i\setminus V_i$ can be shuffled, by~\eqref{e:potatos} we may
assume that
\[
\xi|_{V'_2\setminus V''_2}\quad\text{and}\quad 
\xi|_{V'_3\setminus V''_3}
\]
are isotopic. By~\cite[Lemma~4.13(1)]{GhS} there exists a vertical
convex annulus $A$ with no boundary parallel dividing curve connecting
two ruling curves of $\del V''_2$ and $\del V''_3$. Cutting along $A$
and rounding corners we get a convex vertical torus $T$ surrounding
$F_2$ and $F_3$ with slope $-\frac{1}{p}$. When viewed it as minus the
boundary of the complement of a neighbourhood of $F_1$, the slope of
$T$ becomes $\frac{1}{p}$, which is the critical slope $c_1$.  This
implies that $\xi$ is overtwisted, giving a
contradiction.
\end{proof}

\begin{lem}\label{l:box} 
Let $\xi$ be a positive contact structure on $E_{p,n}$ such that  
\[
\xi|_{\cup_i V'_i} = \xi(q_1,q_2,q_3).
\] 
If $q_1=0$ and $q_3\leq p-1$, or $q_1=1$ and $q_3\geq pn+1$, then
$\xi$ is overtwisted.
\end{lem}

\begin{proof}
We consider the case $q_1=0$ only, because the case $q_1=1$ follows by
a symmetric argument. Assume by contradiction that $\xi$
is tight. Stabilize $F_1$ $n$ times by adding zig-zags to it in such a
way that the newly created basic slices all have negative signs. The
new Legendrian singular fiber has a standard neighbourhood
$V''_1\subset V_1$ such that the boundary slope of $-\del(E_{p,n}\setminus
V''_1)$ is
\[
\frac{n}{pn+1}.
\]
Inside $V_3$ there is a convex neighbourhood $V''_3$ of $F_3$ such that
$-\del(E_{p,n}\setminus V''_3)$ has boundary slope 
\[
-\frac 1{pn+1}.
\]
Moreover, since we can shuffle the basic slices of $V'_3\setminus
V_3$, by the assumption $q_3\leq p-1$ we may assume that 
\[
\xi|_{V'_1\setminus V''_1}\quad\text{and}\quad 
\xi|_{V'_3\setminus V''_3}
\]
decompose into basic slices of the same sign. Therefore,
by~\cite[Lemma~4.13(2)]{GhS} there exists a convex vertical annulus
$A$ between $V''_1$ and $V''_3$ with no boundary parallel dividing
curves. Cutting along $A$ and rounding corners we get a vertical
convex torus which, when viewed as minus the boundary of the
complement of a neighbourhood of $F_2$ has slope $-\frac{n}{pn+1}$,
which is exactly the critical slope $c_2$. This implies that
$\xi$ is overtwisted, giving a contradiction.
\end{proof}

\begin{lem}\label{l:stop} 
Let $\xi$ be a positive contact structure on $E_{p,n}$ such that  
\[
\xi|_{\cup_i V'_i} = \xi(q_1,q_2,q_3).
\] 
If $(q_1,q_2)\in\{(0,0), (1,p)\}$, then $\xi$
is overtwisted for any $q_3\in \{ 0,\ldots, p(n+1)\}$.
\end{lem}

\begin{proof}
Suppose by contradiction that $\xi$ is tight.  Stabilize
$F_1$ $(n+1)$ times and $F_2$ once, and denote by $V''_1$ and $V''_2$ 
standard neighbourhoods of the new Legendrian curves. The slopes of
$-\del(E_{p,n}\setminus V''_1)$ and $-\del(E_{p,n}\setminus V''_2)$ are, 
respectively,
\[
\frac{n+1}{p(n+1)+1}\quad\text{and}\quad
-\frac{n+1}{p(n+1)+1}.
\]
Since $(q_1,q_2)\in\{(0,0),(1,p)\}$, the stabilizations can be
chosen so that 
\[
\xi|_{V'_1\setminus V''_1}\quad\text{and}\quad 
\xi|_{V'_2\setminus V''_2}
\]
decompose into basic slices of the same sign. Therefore,
by~\cite[Lemma~4.13(2)]{GhS} we can find a convex vertical annulus $A$
between $V''_1$ and $V''_2$ with no boundary parallel dividing
curves. Cutting and rounding provides a torus with slope
$\frac{1}{p(n+1)+1}$, which turns into the critical slope $c_3$ when
viewed as minus the boundary of the complement of a neighbourhood of
$F_3$. Therefore, $\xi$ is overtwisted and we have a
contradiction.
\end{proof}

\begin{lem}\label{l:page28}
Let $\xi$ be a positive contact structure on $E_{p,n}$ such that  
\[
\xi|_{\cup_i V'_i} = \xi(q_1,q_2,q_3).
\] 
Suppose that 
\[
(q_1,q_2,q_3)\in\{(0,1,pn+2),(0, p-1, pn+p),
(1,1,0), (1,p-1, p-2)\}.
\]
Then, $\xi$ is overtwisted.
\end{lem}

\begin{proof} 
By contradiction, suppose that $\xi$ is tight. Since the
basic slices of $V'_i\setminus V_i$, $i=2,3$ can be shuffled, the
assumption on $(q_1,q_2,q_3)$ guarantees that we can find convex
neighbourhoods $V''_2$ and $V''_3$ with boundary slope $-\frac{1}{p-1}$
such that $V_i\subset V''_i\subset V'_i$, $i=2,3$, and such that
\[
\xi|_{V'_1\setminus V''_1}\quad\text{and}\quad 
\xi|_{V'_2\setminus V''_2}
\]
are isotopic. Then, by~\cite[Lemma~4.13(1)]{GhS} we can find a convex
vertical annulus between $V''_2$ and $V''_3$ with no boundary parallel
dividing curves. Cutting and rounding gives a convex vertical torus
$T$ which, when viewed as minus the boundary of the complement of a
neighbourhood of $F_1$ has slope $\frac{1}{p-1}$.

Now we follow the line of the argument given in the last paragraph of
the proof of~\cite[Theorem~4.14]{GhS}. By substituting $m_1=1$ into
the formula for the boundary slope $b_1$, we get exactly
$\frac{1}{p-1}$. This shows that $F_1$ can be destabilized to a
Legendrian curve $F'_1$, and $T$ can be viewed as the boundary of a
standard neighbourhood of $F'_1$. If now we stabilize $F'_1$, we get a
new singular fiber $F_1$ and a new standard neighbourhood $V_1$ inside
$V'_1$. But there is a degree of freedom in the choice of the
stabilization of $F'_1$, which corresponds to the choice of
``zig--zag'' to be added to it.  By choosing the appropriate
stabilization, we can arrange a different  sign for the basic slice
$\xi |_{V'_1\setminus V_1}$.

The above argument shows that there is an isotopy between
$\xi$ and a contact structure which restricts to $\cup_i V'_i$ as
$\xi(1-q_1,q'_2,q'_3)$, for some $q'_2$ and $q'_3$ which are apriori
different from $q_2$ and $q_3$. In fact, when we create the torus $T$
we do not touch $V''_2$ and $V''_3$, but we destroy $V'_2\setminus
V''_2$ and $V'_3\setminus V''_3$. Using $-\del(E_{p,n}\setminus V'_1)$,
which has slope infinity, we can find new convex neighbourhoods
$V'_i\supset V''_i$ with infinity boundary slope, but we loose control
on the signs in the basic slice decompositions of $V'_2\setminus
V''_2$ and $V'_3\setminus V''_3$. Since $V''_3$ has been preserved, an
easy computation shows that $q'_3\geq pn+1$ if $q_1=0$, and $q'_3\leq
p-1$ if $q_1=1$. By Lemma~\ref{l:box}, any contact structure which 
restricts to $\cup_i V'_i$ as $\xi(1-q_1,q'_2,q'_3)$ is
overtwisted in these cases and we get a contradiction.
\end{proof}

\begin{proof}[Proof of Theorem~\ref{t:bound}]
Let $V'_i$ $(i=1,2,3)$ be the neighbourhoods given in the statement of
Lemma~\ref{l:twisting}. By~\cite[Lemmas~10, 11]{EH}, there are exactly two 
positive, tight contact structure on $E_{p,n}\setminus\cup_i V'_i$ with
convex boundary and boundary slopes $(\infty, \infty, \infty
)$.  The statement is now an immediate consequence of Lemmas~\ref{l:potatos},
\ref{l:box}, \ref{l:stop} and~\ref{l:page28}.
\end{proof}

\begin{rem}
Shortly after the first version of the present paper was circulated, 
Paolo Ghiggini pointed out to the authors that the upper bound given in 
Theorem~\ref{t:bound} is not sharp for $p>2$.
\end{rem}
 
\section{Generalities in Heegaard Floer theory}
\label{s:three}

In the second part of the paper we will apply Heegaard Floer theory
in proving tightness of certain contact structures (specified by
contact surgery diagrams later) on the oriented 3--manifolds 
$E_{p,n}=S^3 _{p^2n-pn-1}(T_{p,pn+1})$ for $p>1$ and odd.
As it was indicated earlier, the methods used in the subsequent 
sections are completely different from the ones used earlier.
For the sake of completeness
we begin our discussion by shortly reviewing the 
basics of Heegaard Floer theory and contact surgery.

\sh{Ozsv\'ath--Szab\'o homologies}

In a remarkable series of papers \cite{OSzF1, OSzF2, OSzF4, OSz6}
Ozsv\'ath and Szab\'o defined new invariants of many low--dimensional
objects --- including contact structures on closed $3$--manifolds. 
Heegaard Floer theory associates a finetely generated abelian group
$\hf (Y, \t)$ (the \emph{Ozsv\'ath--Szab\'o homology group}) to a
closed, oriented spin$^c$ $3$--manifold $(Y, \t)$, and a homomorphism
\[
F_{W, \s}\co\hf (Y_1, \t_1)\to \hf (Y_2, \t_2)
\]
to an oriented spin$^c$ cobordism $(W, \s)$ between two spin$^c$
$3$--manifolds $(Y_1, \t_1)$ and $(Y_2, \t_2)$. 

Throughout this paper we shall assume that $\Z/2\Z$ coefficients are
being used in the complexes defining the $\widehat{HF}$--groups. With
this assumption, the groups are actually $\Z/2\Z$--vector spaces.  The
group $\hf (Y)$ will denote the sum of $\hf (Y, \t )$ for all spin$^c$
structures. A fundamental property of these groups is that there are
only finitely many spin$^c$ structures on any 3--manifold with
nontrivial Ozsv\'ath--Szab\'o homology groups, hence $\hf (Y)$ is also
finitely generated.  
For a rational homology sphere $Y$ the Ozsv\'ath--Szab\'o homology group
$\hf (Y, \t )$ is nontrivial for any spin$^c$ structure $\t\in Spin ^c (Y)$,
see \cite[Proposition~5.1]{OSzF2}. In particular, for 
a rational homology 3--sphere $Y$ we have  
\[
\dim \hf (Y)\geq \vert H_1(Y; \Z )\vert.  
\]
A rational homology 3--sphere $Y$ is called an~\emph{$L$--space} if
\[
\dim \hf (Y)= \vert H_1(Y; \Z )\vert.
\]
In the light of the above nonvanishing result, this property is equivalent
to 
\[
\hf (Y ,\t )=\Z /2\Z 
\]
for all $\t\in\Spin^c(Y)$.  

Let $Y$ be a closed, oriented 3--manifold and let $K\subset Y$ be a
framed knot with framing $f$. Let $Y(K)$ denote the 3--manifold given
by surgery along $K\subset Y$ with respect to the framing $f$. The
surgery can be viewed at the 4--manifold level as a 2--handle
addition. The resulting cobordism $X$ induces a homomorphism
\[
F_X:=\sum_{\t\in\Spin^c(X)} F_{X,\t}\co \hf (Y)\to\hf(Y(K))
\]
obtained by summing over all spin$^c$ structures on $X$. Similarly,
there is a cobordism $U$ defined by adding a 2--handle to $Y(K)$ along
a normal circle $N$ to $K$ with framing $-1$ with respect to a normal
disk to $K$. The boundary components of $U$ are $Y(K)$ and the
3--manifold $Y'(K)$ obtained from $Y$ by a surgery along $K$ with
framing $f+1$. As before, $U$ induces a homomorphism
\[
F_U\co\hf(Y(K))\to\hf(Y'(K)).
\]
Finally, by attaching a 4--dimensional 2--handle to $Y'(K)$ along a 
normal circle $D$ to $N$ with framing $-1$ with respect to the normal disk
to $N$, we obtain a  cobordism $V$. As it is shown in \cite{LSuj}, the
4--manifold $V$ is a cobordism from $Y'(K)$ to $Y$. As above, 
$F_V$ denotes the induced homomorphism
\[
F_V\co\hf(Y'(K))\to\hf(Y).
\]

\begin{thm} [Surgery exact triangle; \cite{OSzF2}, Theorem~9.16]\label{t:triangle}
The homomorphisms $F_{X}, F_U$ and $F_V$ fit into an exact triangle
\[
\begin{graph}(6,2)
\graphlinecolour{1}\grapharrowtype{2}
\textnode {A}(1,1.5){$\hf (Y)$}
\textnode {B}(5, 1.5){$\hf (Y(K))$}
\textnode {C}(3, 0){$\hf (Y'(K))$}
\diredge {A}{B}[\graphlinecolour{0}]
\diredge {B}{C}[\graphlinecolour{0}]
\diredge {C}{A}[\graphlinecolour{0}]
\freetext (3,1.8){$F_{X}$}
\freetext (4.6,0.6){$F_U$}
\freetext (1.4,0.6){$F_V$}
\end{graph}
\]
\qed
\end{thm}

It was proved in~\cite{OSzF1, OSzabs} that the Ozsv\'ath--Szab\'o
homology groups $\hf (Y)$ split as 
\[
\hf (Y)=\oplus _{(d,\t)\in \J }\hf _d (Y, \t), 
\]
where $\J$ denotes the set of homotopy types of oriented 2--plane
fields on $Y$. The set $\J$ can be identified with $[Y, S^2]$, which
is isomorphic to the set of framed 1--manifolds via the
Pontrjagin--Thom construction. The 1--manifold determines a spin$^c$
structure $\t\in\Spin^c (Y)$, while the framing corresponds to the
degree $d$. This invariant of the oriented 2--plane field $\xi$ is
naturally an element of $\Z / div(\xi )\Z$, where $div(\xi )$ is the
divisibility of $c_1(\xi ) $ in $H^2 (Y; \Z )$.  If $c_1(\xi )$ is
torsion then $div(\xi )=0$. Therefore if $\t \in Spin ^c(Y)$ is
torsion, that is, $c_1(\t )\in H^2 (Y; \Z )$ is a torsion element,
then the Ozsv\'ath--Szab\'o homology group $\hf (Y, \t )$ comes with a
natural relative $\Z$--grading. As it was shown in \cite{OSzabs}, this
relative $\Z$--grading admits a natural lift to an absolute
$\Q$--grading.  In conclusion, for a torsion spin$^c$ structure $\t$ the
Ozsv\'ath--Szab\'o homology group $\hf (Y, \t)$ splits as
\[
\hf (Y, \t )=\oplus _{d \in  \Q}\hf _d (Y, \t), 
\]
where the degree $d$ is determined mod 1 by $\t$.  When $\t
\in\Spin^c(Y)$ has torsion first Chern class, there is an isomorphism
between the homology groups $\hf_d (Y,\t)$ and $\hf _{-d}(-Y,\t)$.

Next we describe the relation between degrees and the maps induced by
4--dimensional cobordisms.  Let $(W, \s )$ be a spin$^c$ cobordism
between two spin$^c$ manifolds $(Y_1,\t_1) $ and $(Y_2, \t _2)$.  If
the spin$^c$ structures $\t_i$ are both torsion and $x\in\hf (Y_1,
\t_1)$ is a homogeneous element of degree $d(x)$, then $F_{W,\s}(x)\in
\hf (Y_2, \t _2)$ is also homogeneous of degree
\[
d(x) + \frac{1}{4}(c_1^2(\s )-3\sigma (W)-2\chi (W)).
\]
Notice that $F_W$ (being equal to the sum $\sum _{\s \in Spin ^c
(W)} F_{W,\s }$) might map a homogeneous element $x\in \hf _d (Y_1, \t
)$ into a nonhomogeneous element $F_W(x) \in \hf (Y_2)$.  

We need one
more piece of information. Recall that the set of spin$^c$ structures
comes equipped with a natural involution, usually denoted by $\t
\mapsto \ot$. The spin$^c$ structure $\ot$, called
the~\emph{conjugate} of $\t$, is defined as follows: If one thinks of
a spin$^c$ structure as a suitable equivalence class of nowhere zero
vector fields (cf. \cite{OSzF1}), then the above involution is the map
induced by multiplying a representative vector field by
$(-1)$. Equivalently, viewing a spin$^c$ structure as an equivalence
class of oriented 2--plane fields, the conjugate action is induced by
reversing the orientation of the planes in the oriented 2--plane field.

\begin{thm}[\cite{OSzF2}, Theorem~2.4]\label{t:conj-iso}
The groups $\hf (Y, \t )$ and $\hf (Y, \ot )$ are canonically isomorphic. \qed
\end{thm}

A spin$^c$ structure $\t \in Spin ^c (Y)$ is induced by a \emph{spin}
structure exactly when $c_1(\t )=0$, or equivalently when $\t = \ot$.
 Let $\J_Y$ denote the isomorphism of Theorem~\ref{t:conj-iso} between
$\hf(Y,\t)$ and $\hf(Y,\ot)$. Then, according
to~\cite[Theorem~3.6]{OSzF4}, given a spin$^c$ cobordism $(W,\s)$ we
have
\begin{equation}\label{e:szam}
F_{W,\s }=\J _{Y'}\circ F_{W, \os }\circ \J _Y,
\end{equation}
where $\overline\s$ is the spin$^c$ structure on the 4--manifold $W$
conjugate to $\s$.  (If we think of $\s\in Spin ^c (W)$ as a suitable
equivalence class of almost--complex structures defined on
$W-\{$finitely many points$\}$, then $\overline\s$ corresponds to the
conjugate of the almost--complex structure defining $\s$.)  As an easy
corollary of (\ref{e:szam}), we get that $F_{W, \s}$ is nontrivial if
and only if $F_{W,\os}$ is nontrivial.  Viewing $\hf (Y)$ with the
conjugate actions as a $\Z / 2\Z$--representation, the above identity
\eqref{e:szam} simply says that the induced map $F_W$ for the
cobordism $W$ is $\Z / 2\Z$--equivariant.

The special relation between spin structures and maps induced by
cobordisms is demonstrated by the following simple observation.
Suppose that $Y$ is a rational homology sphere which is an
$L$--space. We identify the nontrivial element in each group $\hf (Y,
\t )=\Z/2\Z$ with $\t\in Spin ^c (Y)$. With this convention, the set of
spin$^c$ structures provides a basis for $\hf (Y)$. 
Let $V$ be a cobordism between the rational homology spheres
$Y_1$ and $Y_2$, $Y_i$ are $L$--spaces and $\t _i $ are \emph{spin}
structures on $Y_i$ ($i=1,2$).  Let 
\[
{\mathcal {S}}=\{ \s \in\Spin ^c(V) \mid \s \vert_{Y_i} =\t _i \ \  i=1,2\}.
\]
The set ${\mathcal {S}}$ decomposes as the collection ${\mathcal
{S}}_1$ of spin$^c$ structures which are \emph{not} spin structures
and the set of spin structures ${\mathcal {S}}_2$ among the elements
of ${\mathcal {S}}$. As always, let $F_V$ denote the map induced by the
cobordism $V$, that is, $F_V = \sum _{\s \in\Spin ^c (V)}F_{V, \s }$.

\begin{lem}\label{l:ext}
Suppose that $V$ and $\t _i$ ($i=1,2$) are given as above.
If ${\mathcal {S}}_2=\emptyset$ then the $\t _2$--component of
$F_V(\t_1)$ is zero.
\end{lem}

\begin{proof}
Notice that the $\t_2$--component of $F(\t_1)$ is computed by
considering the sum $\sum _{\s \in {\mathcal {S}}}F_{V, \s }(\t _1 )$.
By assumption, this sum is equal to $\sum _{\s \in {\mathcal
{S}}_1}F_{V, \s }(\t _1 )$.  Since ${\mathcal {S}}_1=\{ \s_1,
{\overline {\s_1}}, \ldots , \s_k, {\overline {\s_k}}\} $, $\t
_1={\overline {\t_1}}$ by assumption and $F_{V, \s _i}(\t _1)+ F_{V,
{\overline {\s _i}}}(\t _1)=0$, the lemma follows.
\end{proof}

\sh{Contact $(\pm 1)$--surgery}

Suppose that $L\subset (Y, \xi )$ is a Legendrian knot in a contact
3--manifold. Let $Y_L^{\pm }$ denote the 3--manifold we get by doing
$(\pm 1)$--surgery along $L$, where the surgery coefficient is
measured with respect to the contact framing of $L$. According to the
classification of tight contact structures on a solid torus~\cite{H1},
the contact structure $\xi\vert _{Y-\nu L}$ extends uniquely (up to
isotopy) to the surgered manifolds $Y_L^+$ and $Y_L^-$ as a tight
structure on the glued--up torus. Therefore, the knot $L$ with a
$(+1)$ or $(-1)$ on it uniquely specifies a contact 3--manifold
$(Y_L^+, \xi _L^+)$ or $(Y_L^-, \xi _L^-)$.  (For more about contact
surgery see \cite{DG1, DG2, DGS}.)  In particular, a Legendrian link
${\mathbb {L}} \subset (S^3, \xi _{st})$ in the standard contact
3--sphere (which can be represented by its front projection) defines a
contact structure once the surgery coefficients $(+1)$ and $(-1)$ are
specified on its components. In order to keep diagrams as simple as
possible, we will follow the convention that when in a diagram a
Legendrian knot has no coefficient, then contact $(-1)$--surgery is
performed on it.

\sh{Contact Ozsv\'ath--Szab\'o invariants}

In~\cite{OSz6} Ozsv\'ath and Szab\'o define an invariant
\[
c(Y, \xi)\in \hf (-Y,\t_{\xi })
\]
assigned to a positive, cooriented contact structure $\xi$ on $Y$.  In
fact, $\xi$ (as an oriented 2--plane field) determines an element
$(d(\xi),\t _{\xi })\in \J$ and according to \cite{OSz6} the contact
invariant $c(Y,\xi )$ is an element of $\hf
_{-d(\xi)}(-Y,\t_{\xi})$. Moreover, if $c_1(\xi )\in H^2 (Y; \Z )$ is
torsion then
\[
d(\xi)= \frac{1}{4}(c_1^2(X, J)- 3\sigma (X) -2 \chi (X)+2),
\]
where $X$ is a compact almost--complex 4--manifold with $\partial X=Y$, and
$\xi $ is homotopic to the distribution of complex tangencies on $\partial
X$.

The main properties of the contact Ozsv\'ath--Szab\'o invariant are
summarized in the following two theorems.

\begin{thm}[\cite{OSz6}] \label{t:item1}
If $(Y, \xi )$ is overtwisted, then $c(Y, \xi )=0$. If $(Y, \xi )$ is
Stein fillable then $c(Y, \xi )\neq 0$.  In particular, for the
standard contact structure $(S^3, \xi _{st})$ the invariant $c(S^3, \xi
_{st})\in \hf (S^3)=\Z /2\Z$ is nonzero.\qed
\end{thm}

\begin{thm}[\cite{LS3, OSz6}]\label{t:item2}
Suppose that $(Y_2,\xi_2)$ is obtained from $(Y_1,\xi _1)$ by a
contact $(+1)$--surgery. Then
\[
F_{-W} (c(Y_1, \xi_1))= c(Y_2,\xi_2), 
\]
where $-W$ is the cobordism induced by the surgery with reversed
orientation and $F_{-W}$ is the sum  $\sum_\s  F_{-W,\s}$ over
all spin$^c$ structures on $W$. In particular, if $c(Y_2, \xi_2)\neq
0$ then $(Y_1, \xi_1)$ is tight. \qed
\end{thm}
\noindent Since by~\cite[Proposition~8]{DG1} contact $(-1)$--surgery along a
Legendrian push--off inverts contact $(+1)$--surgery, the above
theorem implies
\begin{cor}\label{c:legsurg}
If $(Y_2,\xi _2)$ is given as Legendrian surgery along a Legendrian knot
in $(Y_1, \xi _1)$ and $c(Y_1, \xi _1)\neq 0$ then $c(Y_2, \xi _2)\neq 0$; 
in particular, $(Y_2, \xi _2)$ is tight. \qed
\end{cor}

\noindent An easy application of the surgery exact triangle and
Theorem~\ref{t:item2} provides
\begin{lem}[\cite{LS3}, Lemma~2.5] \label{l:etak}
The contact structure $\eta _1$ on $S^1\times S^2$ given as contact
$(+1)$--surgery on a Legendrian unknot with Thurston--Bennequin number
$-1$ has nonvanishing contact Ozsv\'ath--Szab\'o invariant
$c(S^1\times S^2, \eta_1)\in \hf (S^1\times S^2)$. \qed
\end{lem}

\section{Symplectic fillings}
\label{s:four}

In this section we show, assuming $n\geq 1$ and $p>1$, that the
3--manifold $E_{p,n}$ does not support fillable contact structures,
thus justifying our use of Heegaard Floer theory in the proof of
tightness of the contact structures described below.

Recall that a compact symplectic 4--manifold $(X, \omega )$ is
a~\emph{symplectic filling} of the closed contact 3--manifold $(Y, \xi
)$ if $\partial X=Y$ and $\omega \vert _{\xi }\neq 0$ along the
boundary $\partial X$.

\begin{prop}\label{p:nofill}
  For each $p>1$ and $n\geq 1$ the oriented 3--manifold
  $E_{p,n}=S^3_{p^2n-pn-1}(T_{p,pn+1})$ is an $L$--space and supports no
  positive, fillable contact structure.
\end{prop}

\begin{proof}
Arguing by contradiction, suppose that $E_{p,n}$ supports a fillable
contact structure. Recall that the slice genus of the $(p,q)$--torus
knot $T_{p,q}$ is equal to $\frac{1}{2}(p-1)(q-1)$.  Since
$(pq-1)$--surgery on the torus knot $T_{p,q}$ is a lens
space~\cite{Mo}, by~\cite[Proposition~4.1]{LSuj} $E_{p,n}$ is an
$L$-space. By~\cite[Theorem~1.4]{OSzUj} this implies that if $(X,
\omega )$ is a symplectic filling of $E_{p,n}$, then $b_2^+(X)=0$. On
the other hand, Figure~\ref{f:newp} shows that $-E_{p,n}$ is the
boundary of a negative definite plumbing 4--manifold $W_{p,n}$.
\begin{figure}[ht]
\begin{center}
\setlength{\unitlength}{1mm}
\begin{picture}(120, 50)
\psfrag{E}{$=$}
\psfrag{R}[l][bl][1][60]{$=$}
\put(0,0){\includegraphics[width=12cm]{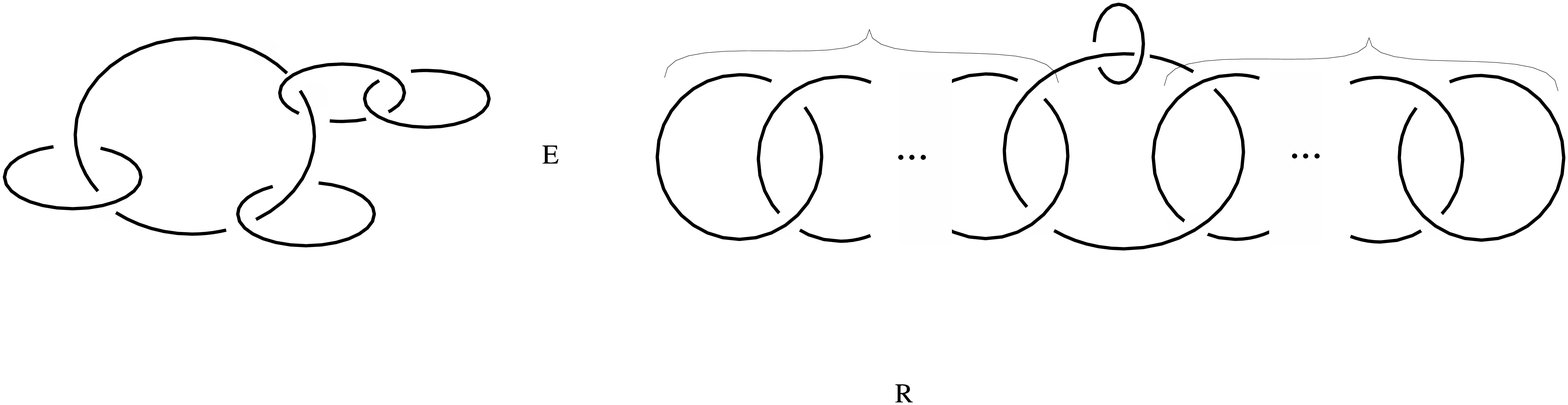}}
\put(5,30){\footnotesize $-p$}
\put(25,27){\footnotesize $p(n+1)+1$}
\put(29,44.5){\footnotesize $p$}
\put(35,43){\ $-n$}
\put(12,43){\footnotesize $0$}
\put(66,47){\footnotesize $p$}
\put(49,28){\footnotesize $-n-1$}
\put(62,28){\footnotesize $-2$}
\put(72,28){\footnotesize $-2$}
\put(82,27){\footnotesize $-2$}
\put(90,28){\footnotesize $-2$}
\put(96,46){\footnotesize $p(n+1)$}
\put(100,28){\footnotesize $-2$}
\put(110,28){\footnotesize $-2$}
\put(85,47){\footnotesize $-p$}
\unitlength=0.8cm
\begin{graph}(15.5,2)(-1.5,-0.4)
\graphnodesize{0.15}
  \roundnode{m3}(-1,0)
  \roundnode{m21}(0.5,0)
  \roundnode{m22}(1.5,0)

  \roundnode{m23}(3.5,0)
  \roundnode{m24}(4.5,0)
  \roundnode{m1}(6.5,0)

  \roundnode{top}(6.5,1)

  \roundnode{n11}(8.5,0)
  \roundnode{n12}(9.5,0)

  \roundnode{n13}(11.5,0)
  \roundnode{n14}(12.5,0)

  \freetext(2.5,0.35){$\overbrace{\hspace{98pt}}$}
  \freetext(2.5,0.7){$p-1$}
  \freetext(10.5,0.35){$\overbrace{\hspace{98pt}}$}
  \freetext(10.5,0.8){$p(n+1)$}

  \edge{m3}{m21}
  \edge{m21}{m22}
  \autonodetext{m22}[e]{\large $\cdots$}

  \autonodetext{m23}[w]{\large $\cdots$}
  \edge{m23}{m24}
  \edge{m24}{m1}
  \edge{m1}{top}
  \edge{m1}{n11}
  \edge{n11}{n12}
  \autonodetext{n12}[e]{\large $\cdots$}

  \autonodetext{n13}[w]{\large $\cdots$}
  \edge{n13}{n14}

  \autonodetext{m3}[s]{\small $-n-1$}
  \autonodetext{m21}[s]{\small $-2$}
  \autonodetext{m22}[s]{\small $-2$}
  \autonodetext{m23}[s]{\small $-2$}
  \autonodetext{m24}[s]{\small $-2$}
  \autonodetext{m1}[s]{\small $-2$}
  \autonodetext{top}[n]{\small $-p$}
  \autonodetext{n11}[s]{\small $-2$}
  \autonodetext{n12}[s]{\small $-2$}
  \autonodetext{n13}[s]{\small $-2$}
  \autonodetext{n14}[s]{\small $-2$}
\end{graph}
\end{picture}    
\end{center}
\caption{Presentation of $-E_{p,n}$ as the boundary of a plumbing}
\label{f:newp}
\end{figure}
Therefore the closed 4--manifold $Z=X\cup _{E_{p,n}}W_{p,n}$ is
negative definite, and by Donaldson's celebrated result~\cite{Do1, Do}
$Z$ has a diagonal intersection form. This implies that any
intersection lattice contained in $Q_{W_{p,n}}$ embeds into the
diagonal intersection form $Q_Z$. But the argument of
\cite[Lemma~4.3]{LS1} with the minor modification given in
\cite[Theorem~4.2]{LSuj} (due to the presence of the framing $-n-1$
instead of $-2$ at the end of one long leg) shows that $Q_{W_{p,n}}$
contains an intersection lattice  which does not embed into
any diagonal intersection form, yielding a contradiction.
\end{proof}

\section{Tight contact structures on $E_{p,n}$}
\label{s:five}

\begin{figure}[ht]
\setlength{\unitlength}{1mm}
\begin{center}
\begin{picture}(100, 80)
%\put(0,0){\tiny\grid(100,80)(5,5)[0,0]}
\psfrag{a}{\footnotesize $p-2$}
\psfrag{b}{\footnotesize $n-1$}
\put(0,39){$\frac{p+3}2$}
\put(60,36){\footnotesize $pn+p$}
\put(63,21){\footnotesize $+1$}
\put(83,40){$\frac{p-1}2$}
\put(0,0){\includegraphics[width=10cm]{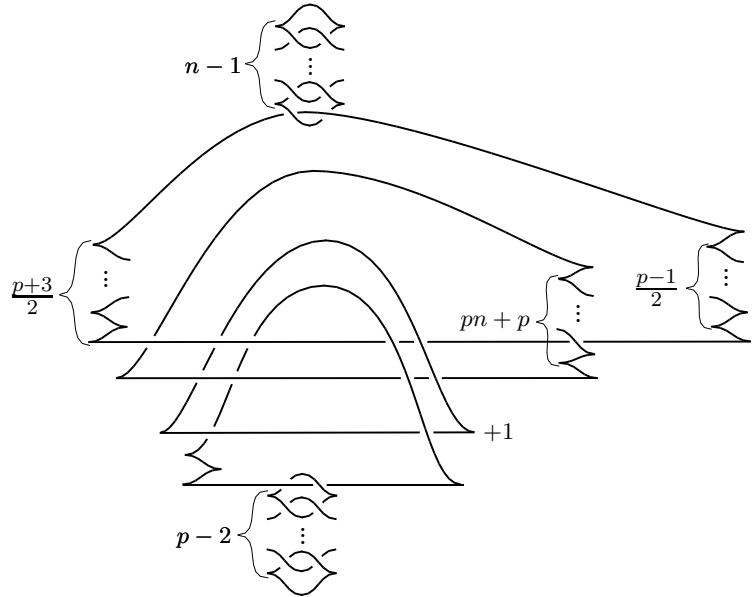}}
\end{picture}
\end{center}
\caption{Tight contact structure on $S_{p,n}$ with $p>1$ odd} 
\label{f:tights}
\end{figure}

Now we outline our approach to the proof of Theorem~\ref{t:T34}.  The
strategy is the following: in this section we specify a contact
structure $\xi _{p,n}$ on a certain 3--manifold $S_{p,n}$ so that the
contact invariant $c(S_{p,n}, \xi _{p,n})$ is nonzero.  Since
$S_{p,n}$ turns out to be an $L$-space, we can identify the invariant
$c(S_{p,n}, \xi _{p,n})\in \hf (-S_{p,n})$ by determining the spin$^c$
structure induced by $\xi _{p,n}$.  By specifying an appropriate
Legendrian knot in $\xi _{p,n}$ and doing contact $(+1)$--surgery
along it, we define a contact structure $\zeta _{p,n}$ on $E_{p,n}$
and a cobordism $X$ from $S_{p,n}$ to $E_{p,n}$.  In the next section
we show that $c(S_{p,n}, \xi _{p,n})$ is not in $\ker F_{-X}$, which
implies that $c(E_{p,n}, \zeta _{p,n})= F_{-X}(c(S_{p,n}, \xi
_{p,n}))$ is nonzero, hence that the contact structure $\zeta _{p,n}$
on $E_{p,n}$ is tight, concluding the argument. Throughout the rest of
the paper we assume that $p>1$ is odd.  The contact structure
$\xi_{p,n}$ is defined by the contact surgery diagram of
Figure~\ref{f:tights}.  The numbers different from $+1$ next to the
vertical braces denote the number of left cusps immediately to their
right.  Moreover (as noted earlier) we adopt the convention that
when in a diagram a Legendrian knot has no coefficient, then contact
$(-1)$--surgery is performed on it.

Notice that the diagram also specifies the underlying oriented
3--manifold $S_{p,n}$.

\begin{prop}\label{p:slides}
The 3--manifold $S_{p,n}$ defined by the contact surgery diagram of 
Figure~\ref{f:tights} is an $L$--space, and the invariant
$c(S_{p,n}, \xi_{p,n})$ is nonzero.
\end{prop}

\begin{proof}
The first statement can be proved in two steps. First observe, by
converting contact surgery coefficients into smooth ones, that
$S_{p,n}$ is orientation preserving diffeomorphic to $S^3_{r}
(T_{p,pn+1})$, with 
\[
r=p(np+1)-\frac{p(n+1)+1}{p(n+1)+2}. 
\]
For the Kirby moves see Figure~\ref{f:traf} and compare the result
with Figure~\ref{f:kirby}.

Since the above $r$ is greater than $2g_s(T_{p,pn+1})-1=p^2n-pn-1$,
by~\cite[Proposition~4.1]{LSuj} the 3--manifold $S_{p,n}$ is an
$L$--space.
\begin{figure}[ht]
\setlength{\unitlength}{1mm}
\begin{center}
\begin{picture}(120, 95)
\psfrag{a}{\footnotesize $a_1$}
\psfrag{e}{\footnotesize $a_2$}
\psfrag{b}{\footnotesize $b$}
\psfrag{c}{\footnotesize $c$}
\psfrag{d}{\footnotesize $d$}
\put(0,0){\includegraphics[width=12cm]{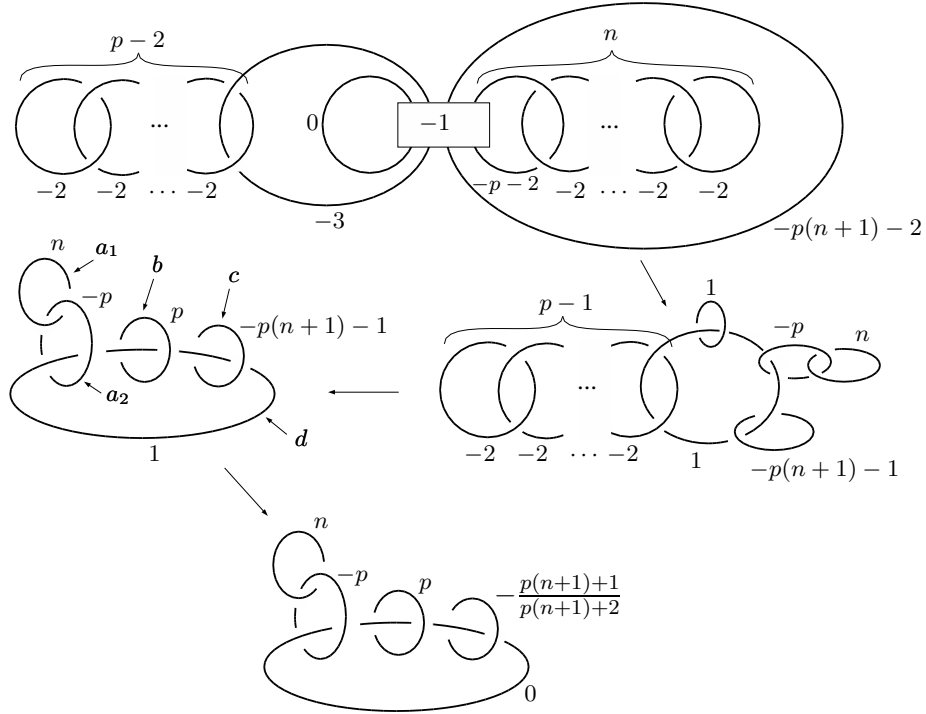}}
%\put(0,0){\tiny\grid(120,95)(5,5)[0,0]}
\put(15,90){\footnotesize $p-2$}
\put(5,70){\footnotesize $-2$}
\put(13,70){\footnotesize $-2$}
\put(20,70){\footnotesize $\cdots$}
\put(25,70){\footnotesize $-2$}
\put(42,66){\footnotesize $-3$}
\put(41,79){\footnotesize $0$}
\put(56,79){\footnotesize $-1$}
\put(63,71){\scriptsize $-p-2$}
\put(74,70){\footnotesize $-2$}
\put(80,70){\footnotesize $\cdots$}
\put(85,70){\footnotesize $-2$}
\put(93,70){\footnotesize $-2$}
\put(80.5,91){\footnotesize $n$}
\put(103,65){\footnotesize $-p(n+1)-2$}
\put(100,33){\footnotesize $-p(n+1)-1$}
\put(103,52){\footnotesize $-p$}
\put(114,51){\footnotesize $n$}
\put(92,34){\footnotesize $1$}
\put(94,57){\footnotesize $1$}
\put(62,35){\footnotesize $-2$}
\put(69,35){\footnotesize $-2$}
\put(76,35){\footnotesize $\cdots$}
\put(81,35){\footnotesize $-2$}
\put(72,55){\footnotesize $p-1$}
\put(32,52){\footnotesize $-p(n+1)-1$}
\put(23,54){\footnotesize $p$}
\put(11,56){\footnotesize $-p$}
\put(7,63){\footnotesize $n$}
\put(20,35){\footnotesize $1$}
\put(45,19){\footnotesize $-p$}
\put(42,26){\footnotesize $n$}
\put(56,18){\footnotesize $p$}
\put(70,3){\footnotesize $0$}
\put(66,16){$-\frac{p(n+1)+1}{p(n+1)+2}$}
\end{picture}
\end{center}
\caption{Surgery diagrams for $S_{p,n}$} 
\label{f:traf}
\end{figure}
The second statement follows from the fact that the structure
$\xi_{p,n}$ is given as Legendrian surgery on the contact structure
$\eta _1$ of Lemma~\ref{l:etak}. Therefore, Lemma~\ref{l:etak} and
Corollary~\ref{c:legsurg} imply that the invariant of $\xi_{p,n}$ is
nonzero.
\end{proof}

\rk{Remark}
In fact, the contact structure $\xi_{p,n}$ can be proved to be Stein
fillable. We will not make use of this fact in our further arguments.

Next, we want to identify the spin$^c$ structure induced by  $\xi_{p,n}$. In
order to do this, we need a little preparation.

It follows from Figure~\ref{f:traf} that the homology group
$H_1(S_{p,n}; \Z )$ has order
\begin{equation}\label{e:hs}
h_S:=|H_1(S_{p,n}; \Z )| = p(pn+1)(p(n+1)+2)-p(n+1)-1.
\end{equation}
Moreover, $H_1(S_{p,n}; \Z )$ is generated by the classes $\mu_{a_1}$,
$\mu_{a_2}$, $\mu_b$, $\mu_c$, $\mu_d$ of suitably oriented meridional
circles to the knots $a_1$, $a_2$, $b$, $c$, $d$ given in
Figure~\ref{f:traf}. These elements are subject to the relations:
\[
n\mu _{a_1} +\mu _{a_2}=0, \quad -p\mu _{a_2}+\mu _{a_1}+\mu _d=0, \quad
p\mu _b +\mu _d=0, 
\]
\[
(-p(n+1)-1)\mu _c +\mu _d=0, \quad 
\mu _{a_2}+\mu _b +\mu _c +\mu _d=0.
\]
The relations above imply that $\mu_d$ generates the homology group, since
$\mu_{a_1}$, $\mu_{a_2}$, $\mu_b$ and $\mu_c$ can be expressed in terms of
$\mu_d$ as
\begin{itemize}
\item $\mu_{a_1}=[n(n+1)p^2+2np-1-n]\mu _d$, $\ \ \ \ \mu _{a_2}=-n\mu_{a_1}$,
\item $\mu_b =[(-n^2-n)p^2+(-1-3n)p-1+n]\mu _d$,
\item $\mu_c=[(n^2+2n+1)np^2+p(2n^2+3n+1)-(n+2)n]\mu _d$.
\end{itemize}

Notice that the order of $H_1 (S_{p,n}; \Z)$ is always odd. Therefore, there
is no 2--torsion in the second cohomology of $S_{p,n}$, and the spin$^c$
structures on $S_{p,n}$ are determined by their first Chern classes.

\begin{lem}\label{l:spinccomp}
Let $\t_{p,n}$ be the spin$^c$ structure induced by $\xi_{p,n}$. Then,
if $p$ is odd we have $c_1(\t_{p,n})=PD(\mu_d)$.
\end{lem}

\begin{proof}
Consider the 4--manifold $X$ determined by the surgery diagram of
Figure~\ref{f:tights}. Since $X$ is simply connected, a spin$^c$
structure on $X$ is determined by its first Chern class. Let $\al\in
H^2(X;\Z)$ be the unique cohomology class which evaluates on each
2--homology class corresponding to an oriented knot $K$ of the diagram
as the rotation number of $K$. Then, the spin$^c$ structure
corresponding to $\al$ restricts to the spin$^c$ structure of
$\xi_{p,n}$ (see e.g.~\cite{DGS} for details).

Therefore, after choosing a suitable orientation of the curves in
Figure~\ref{f:tights}, we have
\begin{equation}\label{e:pddual}
PD(c_1(\t_{p,n })) = \sum_K \rot(K) \mu_K,
\end{equation}
where the sum is over all surgery curves, $\rot(K)$ denotes the
rotation number of the oriented Legendrian knot $K$ and $\mu_K$
denotes the first homology class induced by its meridian. Recall that
according to \cite{Go, GS} the front projection determines the
rotation number of the corresponding Legendrian knot as
\begin{equation}\label{e:rot}
{\rm {rot}}(K)=\frac 12(c_d -c_u),
\end{equation}
where $c_u $ and $c_d$ denote the number of up and down cusps in the
projection. Using Formulas~\eqref{e:pddual} and~\eqref{e:rot}, and
following the Kirby moves of Figure~\ref{f:traf}, one can
easily check that
\[
PD(c_1(\t_{p,n }))= -\mu_{a_2}-\mu_b +p(n+1)\mu_c-\mu_d.
\]
Replacing each of $\mu_{a_2}$, $\mu_b$ and $\mu_c$ by the
corresponding multiple of $\mu_d$ yields, after a somewhat tedious
calculation, $PD(c_1(\t_{p,n }))=\mu_d$.
\end{proof}

\begin{figure}[ht]
\setlength{\unitlength}{1mm}
\begin{center}
\begin{picture}(130, 110)
\psfrag{a}{\footnotesize $n$}
\psfrag{b}{\scriptsize $-p-2$}
\psfrag{c}{\footnotesize $p$}
\psfrag{d}{\footnotesize $-p(n+1)-1$}
\psfrag{0}{\footnotesize $0$}
\psfrag{e}{\footnotesize $p-2$}
\psfrag{f}{\footnotesize $-2$}
\psfrag{g}{\footnotesize $\cdots$}
\psfrag{h}{\footnotesize $-3$}
\psfrag{i}{\footnotesize $-1$}
\psfrag{m}{\footnotesize $n$}
\psfrag{n}{\footnotesize $-p(n+1)-2$}
\psfrag{p}{\footnotesize $-p$}
\psfrag{L}{\footnotesize $L$}
\psfrag{z}{$\frac{p-1}2$}
\psfrag{w}{\footnotesize $n-1$}
\put(0,0){\includegraphics[width=13cm]{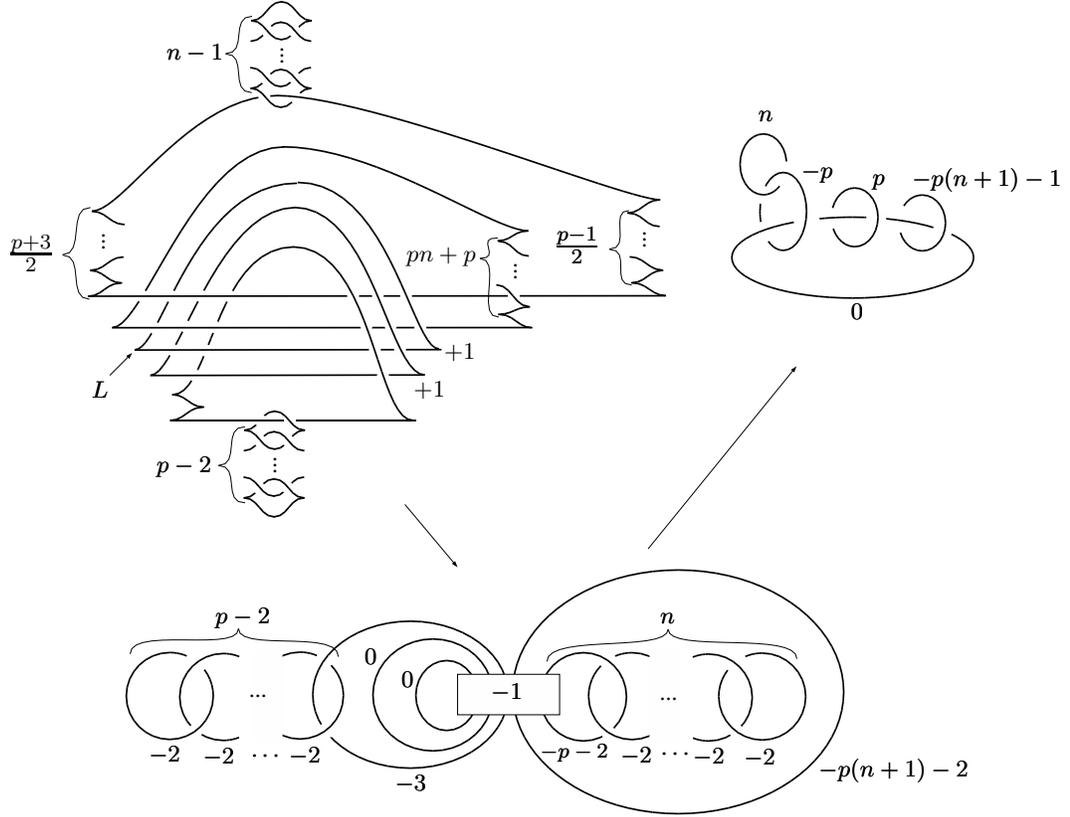}}
\put(1,74){$\frac{p+3}2$}
\put(54,74){\footnotesize $pn+p$}
\put(55,56){\footnotesize $+1$}
\put(59,61){\footnotesize $+1$}
\end{picture}
\end{center}
\caption{The contact structure $\zeta_{p,n}$ on $E_{p,n}$}
\label{f:legknot}
\end{figure}

\begin{defn}
Let $\zeta_{p,n}$ be the contact structure defined by the upper--left 
contact surgery picture of Figure~\ref{f:legknot}.
\end{defn}

\begin{prop}
The contact structure $\ze_{p,n}$ is supported by $E_{p,n}$.
\end{prop}

\begin{proof}
The proof requires only a minor modification of the Kirby calculus of
Figure~\ref{f:traf}. This modification is shown in
Figure~\ref{f:legknot}.
\end{proof}

\section{Maps between the Ozsv\'ath--Szab\'o homologies}
\label{s:six}

\begin{figure}[ht]
\begin{center}
\psfrag{a}{\footnotesize $-1$}
\psfrag{b}{\footnotesize $-n$}
\psfrag{c}{\footnotesize $p$}
\psfrag{d}{\footnotesize $-p$}
\psfrag{0}{\footnotesize $0$}
\psfrag{e}{\footnotesize $p(n+1)+1$}
\psfrag{f}{\footnotesize $a_1$}
\psfrag{g}{\footnotesize $d$}
\psfrag{h}{\footnotesize $a_2$}
\psfrag{l}{\footnotesize $b$}
\psfrag{m}{\footnotesize $c$}
\psfrag{X}{\footnotesize $-X$}
\psfrag{V}{\footnotesize $V$}
\psfrag{S}{\footnotesize $-S_{p,n}=$}
\psfrag{E}{\footnotesize $=-E_{p,n}$}
\psfrag{L}{\footnotesize $=-L_{p,n}$}
\psfrag{K}{\footnotesize $K$}
\includegraphics[width=12cm]{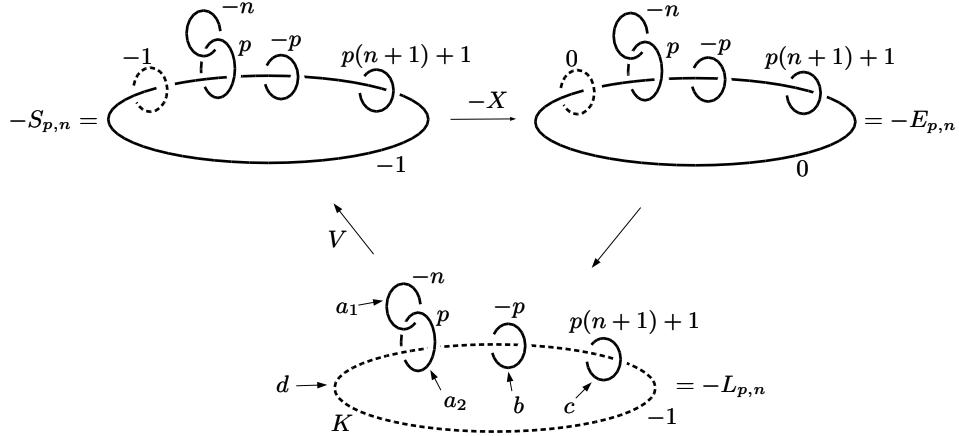}
\end{center}
\caption{Manifolds and cobordisms in the main surgery triangle}
\label{f:tri}
\end{figure}

In this section we show that the contact Ozsv\'ath--Szab\'o invariant
$c(E_{p,n}, \zeta _{p,n})$ is nonzero. This proves
Theorem~\ref{t:T34}.  Note that $\zeta _{p,n}$ is obtained by contact
$(+1)$--surgery on $\xi_{p,n}$ along the Legendrian knot $L$ shown in
Figure~\ref{f:legknot}. There is a cobordism naturally associated to
the surgery which we denote by $X$. By the properties of the contact
Ozsv\'ath--Szab\'o invariants we know that $c(E_{p,n}, \zeta _{p,n})=
F_{-X}(c(S_{p,n}, \xi _{p,n}))$. This section is devoted to collect
partial information about the map $F_{-X}$. In particular, we show
that $c(S_{p,n}, \xi_{p,n})$ is not in $\ker F_{-X}$. Recall that we
have assumed that $p>1$ is odd.  The cobordism $-X$ induced by the
surgery on the knot $L$ of Figure~\ref{f:legknot} (after reversing its
orientation) fits into the triangle given by Figure~\ref{f:tri}.  
In the remaining figures of the paper we adopt the convention of
denoting the 3--manifold under examination by solid framed links,
while dashed curves denote the 2--handles of the cobordism built on
the given 3--manifold. We shall use the corresponding exact triangle
involving the Ozsv\'ath--Szab\'o homology groups to study the map
\[
F:=F_{-X}\colon \hf (-S_{p,n})\to \hf (-E_{p,n}).
\]
The strategy to show that the contact invariant 
\[
c(E_{p,n}, \zeta_{p,n})=F_{-X}(c(S_{p,n}, \xi _{p,n})) 
\]
is nonzero will be the following. Let $G_V$ be the map induced by the
cobordism $V$. First we show that there exists an element of
$\hf(-L_{p,n})$ corresponding to a spin structure on $-L_{p,n}$ with
the property that its $G_V$--image is equal to $a+{\overline {a}}$ for
some $a\in \hf (-S_{p,n})$.  (Recall that ${\overline {a}}$ denotes
the image of $a\in \hf (-S_{p,n})$ under the ${\mathcal {J}}$--action
induced by conjugation on spin$^c$ structures.)  Next we consider the
decomposition of this element $a$ into a sum of homogeneous terms, and
we find a homogeneous component $a_1\in \hf (-S_{p,n}, \t )$ which
maps to a nonzero element under $F_{-X}$.  In the final step of the
proof we determine the spin$^c$ structure $\t$ corresponding to the
above element $a_1$ and show that it is equal to the spin$^c$
structure induced by the contact structure $\xi _{p,n}$. Since
$S_{p,n}$ was proved to be an $L$--space, the nonzero elements $a_1$
and $c(S_{p,n}, \xi _{p,n})$ inducing the same spin$^c$ structure must
be equal. In particular, $F_{-X}(c(S_{p,n}, \xi _{p,n}))\neq 0$,
concluding the proof.  In identifying the spin$^c$ structure of the
element $a_1$ we appeal to a computation which determines the degree
difference between two spin structures on $-L_{p,n}$ and $-E_{p,n}$;
this computation relies on the study of a related exact triangle and
is given in a separate subsection.  Notice that all the 3--manifolds
in the triangle of Figure~\ref{f:tri} are $L$--spaces: this property
was verified for $E_{p,n}$ and $S_{p,n}$ in
Propositions~\ref{p:nofill} and \ref{p:slides}, while $L_{p,n}$ is the
connected sum of three lens spaces, hence the $L$--space property
trivially follows.  (Recall that $\hf (Y)$ is isomorphic to $\hf (-Y)$
hence $Y$ is an $L$--space if and only if $-Y$ is an $L$--space.)  To
set up notation, consider the surgery exact triangle defined by the
cobordisms of Figure~\ref{f:tri}:
\begin{equation}\label{e:triangle}
\begin{graph}(6,2)
\graphlinecolour{1}\grapharrowtype{2}
\textnode {A}(1,1.5){$\hf (-S_{p,n})$}
\textnode {B}(5, 1.5){$\hf (-E_{p,n})$}
\textnode {C}(3, 0){$\hf (-L_{p,n})$}
\diredge {A}{B}[\graphlinecolour{0}]
\diredge {B}{C}[\graphlinecolour{0}]
\diredge {C}{A}[\graphlinecolour{0}]
\freetext (3,1.8){$F=F_{-X}$}
\freetext (4.4,0.6){$H$}
\freetext (1.5,0.6){$G_V$}
\end{graph}
\end{equation}
Using the surgery descriptions it follows that
\begin{equation}\label{e:he}
h_E:=\vert H_1(E_{p,n}; \Z )\vert = p^2n-pn-1, \mbox{ and} 
\end{equation}
\begin{equation}\label{e:hl}
h_L:=\vert H_1(L_{p,n})\vert = p(pn+1)(p(n+1)+1).
\end{equation}

\begin{prop}\label{p:van}
The map $H$ is equal to 0, therefore $F$ is surjective and $G_V$ is 
injective.
\end{prop}

\begin{proof}
Since the three 3--manifolds are all $L$--spaces, their
Ozsv\'ath--Szab\'o homology groups can be determined from their first
homologies. Now a simple computation using
Equations~\eqref{e:hs},~\eqref{e:he} and~\eqref{e:hl} shows that
$h_E+h_L=h_S$, hence the statement of the lemma follows from the
exactness of the triangle and elementary algebra (cf.~also the
concluding remark of~\cite[Section~2]{LSuj}).
\end{proof}

\begin{lem}
The manifolds $S_{p,n}$ and $E_{p,n}$ admit a unique spin structure, 
while $L_{p,n}$ supports exactly two spin structures.
\end{lem}

\begin{proof}
Recall that any orientable 3--manifold $Y$ admits a spin structure, and
the number of inequivalent spin structures is given by $\vert H^1 (Y;
\Z /2\Z )\vert$. Using Equations~\eqref{e:hs},~\eqref{e:he}
and~\eqref{e:hl} it is easy to check that $S_{p,n}$ and $E_{p,n}$
have first homology groups of odd order, while for $L_{p,n}$ 
(as the connected sum of the three lens spaces of Figure~\ref{f:tri}) we have 
$H^1(L_{p,n}; \Z/2\Z )=\Z /2\Z$.
\end{proof}

\begin{lem}\label{l:spinonL}
Let $V$ and $W$ be the cobordisms defined, respectively in
Figure~\ref{f:tri} and Figure~\ref{f:reltri}. Then, each spin
structure on $-L_{p,n}$ extends as a spin structure to one of the
cobordisms $V$ and $W$, but not to the other.
\end{lem}

\begin{proof}
Recall that we are assuming that $p$ is odd. In the proof we will
distinguish two cases according to the parity of $n$. We would like to
present $-L_{p,n}$ as the boundary of two spin 4--manifolds.  Consider
the bottom pictures of Figures~\ref{f:tri} and~\ref{f:reltri}. Suppose
first that $n$ is even. By anti-blowups we can transform the
$(-n)$--framed unknot linking the $p$--framed unknot into a chain of
$(+2)$'s. During this operation we change the framing $p$ into
$p+1$. Do the same operation with the $(-p)$--framed circle. Notice
that after the above blow ups and blow downs the parity of the framing
of the knot $K$ shown by the figures has changed.  Since $n$ is even,
$p(n+1)+1$ is also even.  Therefore the diagram defines a simply
connected spin 4--manifold with a unique spin structure, and we define
$\t_V\in\Spin (L_{p,n})$ as the restriction of this unique spin
structure to the boundary.  Since the framing of $K$ when defining $V$
is even, $\t_V$ extends to $V$ as a spin structure but does not extend
to $W$ (as a spin structure), since it would give a spin 4--manifold
with a homology class of odd square, hence with nontrivial second
Stiefel--Whitney class.

To find the other spin 4--manifold, we turn the $(p(n+1)+1)$--framed
circle into a chain of $(-2)$'s by blowing up and down. This operation
changes the parity of the framing of $K$ again. We define $\t_W$ as
the restriction of the unique spin structure of the resulting simply
connected spin 4--manifold. Since the parity of the framing of $K$ is
now different than in the previous case, the spin structure $\t_W$
extends to the cobordism $W$ as a spin structure but does not extend
to $V$ as a spin structure. Clearly $\t_V\neq\t_W$, and when $n$ is
even we are done.

Finally we address the case of odd $n$. In this case both $-p$ and
$p(n+1)+1$ are odd, so first we turn these surgeries into chains of
$(+2)$ (and $(-2)$, resp.) surgeries. Each one of these
transformations changes the framing of the knot $K$ by $+1$ (and $-1$
resp.), so the net change of the framing of $K$ is zero. Now we have a
choice for the remaining two odd framed surgery curve defining
$-L_{p,n}$.  If we turn the $(-n)$--framed unknot into a chain of
$(+2)$'s, we change the framing $p$ into $p+1$, but we do not change
the framing of $K$.  Hence the resulting 4--manifold admits a spin
structure  $\s_W$ which extends to $W$ as a spin structure, but
not to $V$. We denote the restriction of $\s_W$ to the boundary
$-L_{p,n}$ by $\t_W$. On the other hand, the corresponding operation
on the $p$--framed circle changes the framing of the $(-n)$--framed
circle to $(-n-1)$ and also changes the parity of the framing of $K$.
Therefore the spin structure of the resulting simply connected spin
4--manifold will extend to $V$ as a spin structure but not to $W$.
The restriction of this spin structure to the boundary $-L_{p,n}$ will
be called $\t_V$. Clearly $\t_W\neq\t_V$, and the proof is finished.
\end{proof}

%\rk{Remark.} The same statement, and a very similar argument apply
%when $p$ is even, but we will not make any use of this parallel
%statement.

\rk{Notation.}  We denote the unique spin structures on $-S_{p,n}$ and
$-E_{p,n}$, respectively, by $\t_S$ and $\t_E$. 
As in the proof of Lemma~\ref{l:spinonL}, we denote by $\t_V$
the spin structure on $-L_{p,n}$ which extends as a spin structure to
$V$ but not to $W$ and by $\t_W$ the spin structure which extends 
(as a spin structure) to $W$ but not to $V$.

\subsection*{Computations}

Now we return to the analysis of Triangle~\eqref{e:triangle}. 
Recall that when $Y$ is a rational homology sphere which is an
$L$--space, we have identified the nontrivial element in each group $\hf (Y,
\t )$ with $\t\in Spin ^c (Y)$. If $H_1(Y; \Z )$ is of odd rank, then
$Y$ admits a unique spin structure, which will be denoted by
$\t_Y$. Using the conjugate action encountered in
Section~\ref{s:three} (cf. Theorem~\ref{t:conj-iso}), 
and denoting ${\mathcal {J}}(\t ) $ by ${\overline {\t}}$, 
in this case the
vector space $\hf(Y)$ has a basis of the form
\begin{equation}\label{e:basis}
\{\t_1,\overline{\t_1},\t_2,\overline{\t_2},\ldots,\t_k,\overline{\t_k},\t_Y\}.
\end{equation}
Let
\[
C :=\langle \t_1,\ldots, \t_k\rangle\subset\hf(Y).
\]
Then, we have 
\begin{equation}\label{e:dsum}
\hf(Y)=\langle\t_Y\rangle\oplus C\oplus\overline{C}.
\end{equation}
Notice that the subspace $C\subset \hf (Y)$ depends on a choice of
basis as in~\eqref{e:basis}, therefore the above splitting is not
canonical. In analogy to
Equation~\eqref{e:dsum}, there are direct sum decompositions
\begin{eqnarray}\label{e:decomp}
\hf(-S_{p,n}) = \langle \t_S\rangle\oplus A\oplus\overline{A},\notag\\
\hf(-L_{p,n})= \langle \t_V\rangle\oplus\langle \t_W\rangle\oplus
C\oplus\overline{C}\\
\hf(-E_{p,n})=\langle \t_E \rangle\oplus T\oplus\overline{T}.\notag
\end{eqnarray}
Since by its definition $\t_W$ does not extend as a spin structure to
$V$, Lemma~\ref{l:ext} implies that
\[
G_V(\t_W)\in A\oplus\overline{A}. 
\]
Since $\t_W$ is fixed under conjugation, so is $G_V(\t_W)$, therefore 
there is an element $a\in A$ such that $G_V(\t_W)=a+\oa $. 
Notice that $F(a)=F(\oa)=\overline{F(a)}$, because
\[
F(a)+F(\oa) = F(a+\oa ) = F(G_V(\t_W)) = 0,
\]
and we work with $\Z / 2\Z$--coefficients.
\begin{lem}
We have $F(a)\neq 0$.
\end{lem}

\begin{proof}
If $F(a)=0$, then by exactness
$a=G_V(c)$ for some $c\in\hf(-L_{p,n})$. Therefore 
\[
G_V(\t_W+c+\overline{c})=0.
\]
Since $c+\overline{c}\in C\oplus\overline{C}$, the injectivity of
$G_V$ would imply $\t_W\in C\oplus\overline{C}$, which is impossible
by~\eqref{e:decomp}.
\end{proof}

\begin{lem}\label{l:tE-comp}
Suppose that $F(a)=\epsilon \t_E+t +\overline {t}$ for some
$t\in \hf (-E_{p,n})$. Then, $\epsilon \neq 0$.
\end{lem}

\begin{proof}
By contradiction, suppose that $\epsilon =0$. By the surjectivity of $F$,
there is $b\in \hf (-S_{p,n})$ with $F(b)=t$, implying also
$F({\overline {b}})={\overline {t}}$. Now consider $x=a+b +{\overline
{b}}$. Then, $F(x)=0$, and so $F({\overline {x}})=0$. By exactness
this means that there is $u\in \hf (-L_{p,n})$ satisfying $G_V(u)=x$,
and so $G_V ({\overline {u}})={\overline {x}}$.  This implies that
$G_V(u+{\overline {u}}+\t_W)=0$. By the the injectivity of $G_V$, this
would imply
\[
\t_W=u+\overline{u}\in C\oplus\overline{C},
\]
which is impossible by~\eqref{e:decomp}.  
\end{proof}

In order to apply the degree--shift formula for the cobordisms $X$ and $V$,
we need some understanding of their algebraic topology.

\begin{lem}\label{l:signatures}
We have 
\[
H_2(V;\Z)\cong H_2(-X;\Z)\cong\Z
\]
and 
\[
\si(V)=\si(-X)=-1,
\]
where $\si$ denotes the signature.
\end{lem}

\begin{proof}
The cobordism $V$ is obtained by attaching a 2--handle to the rational
homology sphere $-L_{p,n}$. Therefore, $H_2(V,-L_{p,n};\Z)\cong\Z$, and 
the exactness of the sequence 
\[
0\lra H_2(V;\Z)\lra H_2(V,-L_{p,n};\Z)\lra H_1(-L_{p,n};\Z)
\]
implies that $H_2(V;\Z)\cong\Z$. A similar argument shows that
$H_2(-X;\Z)\cong\Z$.

It is easy to deduce from Figure~\ref{f:tri} that 
\[
V\cup -X\cong Q\#\overline{\CP^2},
\]
where the cobordism $Q\#\overline{\CP^2}$ is given by
Figure~\ref{f:Y}, obtained by applying two Rolfsen twists to the
bottom picture of Figure~\ref{f:tri}.
\begin{figure}[ht]
\begin{center}
\psfrag{a}{\footnotesize $-1$}
\psfrag{b}{\small $-\frac{p(n+1)+1}{p(n+1)}$}
\psfrag{c}{\footnotesize $-p$}
\psfrag{d}{\small $-\frac{np+1}{n(p-1)+1}$}
\psfrag{C}{\footnotesize $\cong$}
\psfrag{e}{\footnotesize $-2$}
\psfrag{f}{\footnotesize $-3$}
\includegraphics[width=12cm]{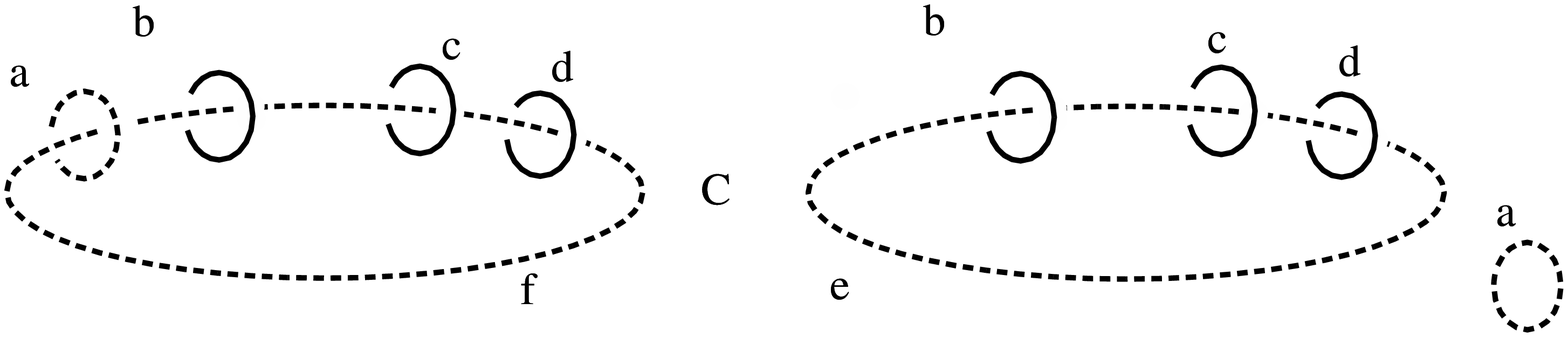}
\end{center}
\caption{The cobordism $Q\#\overline{\CP^2}$}
\label{f:Y}
\end{figure}
As in the proof of Lemma~\ref{l:spinonL} we can replace the two
unknots with non--integral surgery coefficients by two chains of
unknots with integral coefficients, with each coefficient less then or
equal to $-2$. The resulting picture expresses $Q$ as a 4--dimensional
2--handle attached to the boundary of a 4--dimensional plumbing $P$
with $\del P=-L_{p,n}$. Moreover, the union $P\cup Q$ is still a
plumbing and we claim that it is negative definite. In fact, according
to~\cite[Theorem 5.2]{NR}, to see this it is enough to check that
\[
-2 + \frac{n(p-1)+1}{np+1} + \frac 1p + \frac{p{n+1}}{p(n+1)+1} < 0
\]
for any $n\geq 1$ and $p\geq 2$. This implies that $Q$ is negative definite
and concludes the proof.
\end{proof}

Recall that $h_S, h_E$ and $h_L$ denote the cardinality of the
homology groups $H_1(S_{p,n}; \Z)$, $H_1(E_{p,n}; \Z )$ and
$H_1(L_{p,n};\Z )$, respectively. 

\begin{lem}\label{l:square}
Let $g\in H_2(V;Z)$ and $g'\in H_2(-X;\Z)$ be generators. Then, 
\[
g\cdot g = - h_L h_S,\quad\text{and}\quad g'\cdot g' = -h_S h_E.
\]
\end{lem}

\begin{proof}
We give the argument for $V$, the one for $-X$ being essentially the
same. From Figure~\ref{f:tri} we see that $V$ is obtained by attaching
a 4--dimensional 2--handle along a circle which represents a generator
of $H_1(-L_{p,n};\Z)$. Therefore, $H_1(V;\Z)=0$. Since by
Lemma~\ref{l:signatures} we have that $H_2(V;\Z)\cong\Z$, the
universal coefficient theorem gives
\[
H_2(V,\del V;\Z)\cong H^2(V;\Z)\cong\Z.
\]
Consider the exact sequence 
\begin{equation}\label{e:pairseq}
0\to H_2(V;\Z)\stackrel{i_*}{\to} H_2(V,\del V;\Z)\to H_1(\del V;Z)
\cong\Z/h_L\Z\oplus \Z/h_S\Z\to 0.
\end{equation}
It is easy to check that $h_L$ and $h_S$ are coprime, thus 
\[
\Z/h_L\Z\oplus \Z/h_S\Z\cong \Z/(h_L h_S)\Z,
\]
and $i_*(g)$ must be equal to $h_L h_S$ times a generator of
$H_2(V,\del V;\Z)$. Therefore, since by Lemma~\ref{l:signatures} 
the cobordism $V$
has negative definite intersection form,
\[
g\cdot g = \langle \PD(i_*(g)), g\rangle = -h_L h_S.
\]
\end{proof}

\begin{lem}\label{l:c1V}
Let $\s\in\Spin^c(V)$, and let $C\subset V$ be the cocore of 
the 2--handle defining $V$. If $\s|_{-L_{p,n}}=\t_W$, then 
\[
\PD(c_1(\s)) = k [C] \in H_2(V,\del V;\Z)
\]
for some odd integer $k$. Moreover, 
\[
c_1(\s)\cdot c_1(\s) = -\frac{k^2 h_L}{h_S}.
\]
\end{lem}

\begin{proof}
According to the proof of Lemma~\ref{l:spinonL}, $\t_W$ is the
restriction to $-L_{p,n}$ of a spin structure $\u$ on a spin 4--manifold
$Z$ with $\del Z=-L_{p,n}$.  Moreover, $Z$ is obtained by attaching
4--dimensional 2--handles to the 4--ball $B^4$, and $V$ by attaching a
last 2--handle $H$ to $\del Z$. Recall that the framing of the
attaching circle of $H$ is odd, because $\t_W$ does not extend over $V$ as a
spin structure. Thus, if $\s|_{-L_{p,n}}=\t_W$, then $\s$ extends $\u$
to $W:=Z\cup V$ as a spin$^c$ structure. Denote by $\tilde\s$ the
extended spin$^c$ structure $\u\cup\s$. Thinking of $H$ as attached to
$S^3=\del B^4$, let $F$ denote the surface obtained by capping off the
core $D$ of $H$ by a Seifert surface with interior pushed in $B^4$. Since
$c_1(\tilde\s)$ is characteristic and $F$ has odd square, we have
\[
\langle c_1(\tilde\s), [F]\rangle = k 
\]
for some odd integer $k$. Therefore, since $W$ is simply connected,
$\PD(c_1(\tilde\s))=k[C]$. The first part of the statement follows
because $\tilde\s$ restricts to $\s$ on $V$ and $C\subset V$.

Now observe that the boundary of $h_L$ parallel copies of $D$ is homologically
trivial in $-L_{p,n}$. Thus, we can define $S\subset V$ to be the surface
obtained by capping off $h_L D$ in $-L_{p,n}$ with a bounding surface.
Moreover, since $C$ is disjoint from $-L_{p,n}$, by Exact
Sequence~\eqref{e:pairseq} the relative homology class $[C]$ must be a
multiple of $h_L$ times a generator $g'$ of $H_2(V,\del V;\Z)$. But the
equality $[C]\cdot [S]= h_L$ implies at once that $[S]$ is a generator $g$ of
$H_2(V;\Z)$, and $[C]$ is $h_L g'$. Now recall that in the proof of
Lemma~\ref{l:square} we showed that the image of $g$ under the map $i_*$ of
Exact Sequence~\eqref{e:pairseq} is equal to $\pm h_L h_S g'$.  Therefore,
\[
h_S\PD(c_1(\s)) = k h_S [C] = k h_S h_L g' = \pm k i_*(g)
\]
which implies, by Lemma~\ref{l:square}, that
\[
c_1(\s)\cdot c_1(\s) = k^2\frac {g\cdot g}{h_S^2} 
= -k^2\frac {h_L}{h_S}.
\]
\end{proof}

\begin{lem}\label{l:c1X}
Let $\s\in\Spin^c(-X)$, and let $D\subset -X$ be the core of 
the 2--handle defining $-X$. If $\s|_{-E_{p,n}}=\t_E$, then 
\[
\PD(c_1(\s)) = l [D] \in H_2(-X,\del (-X);\Z)
\]
for some odd integer $l$. Moreover, 
\[
c_1(\s)\cdot c_1(\s) = -l^2 \frac{h_E}{h_S}.
\]
\end{lem}

\begin{proof}
  Observe that the spin structure $\t _E$ does not extend to $-X$
  simply because $-X$ does not carry spin structures. This follows
  immediately from Lemma~\ref{l:square}, since both $h_S$ and $h_E$
  are odd numbers. Thus, the proof of this lemma is similar to the
  proof of Lemma~\ref{l:c1V}, and we omit it.
\end{proof}

We wish to find a relation between the degrees of $\t_W$ and $\t_E$.
This can be done with a (quite tedious) direct computation: the
gradings of generators of $\hf (Y)$ for a lens space $Y$ are given
in~\cite{OSzabs}, and since $-L_{p,n}$ is a connected sum of three
lens spaces and the degrees are additive under connected sums, the
computation of the degree of $\t_W$ is a fairly easy exercise.  The
degree of an element in the Ozsv\'ath--Szab\'o homology of a Seifert
fibered 3--manifold can be computed using formulae from~\cite{nem,
OSzplum}. In particular, in~\cite{nem} there is an explicit formula in
terms of a vector with some special properties in the cohomology of a
certain negative definite plumbing with boundary $Y$. This direct
computation, however, is quite delicate, so we prefer to choose a
theoretically more involved, less computational way of relating the
degrees of $\t_W$ and $\t_E$. In particular, we will get the desired
conclusion by studying a related triangle of manifolds.

\sh{Digression: study of a related triangle}

\begin{figure}[ht]
\begin{center}
\psfrag{a}{\footnotesize $-1$}
\psfrag{b}{\footnotesize $-n$}
\psfrag{c}{\footnotesize $p$}
\psfrag{d}{\footnotesize $-p$}
\psfrag{0}{\footnotesize $0$}
\psfrag{1}{\footnotesize $1$}
\psfrag{e}{\footnotesize $p(n+1)+1$}
\psfrag{W}{\footnotesize $W$}
\psfrag{U}{\footnotesize $=-U_{p,n}$}
\psfrag{E}{\footnotesize $-E_{p,n}=$}
\psfrag{L}{\footnotesize $=-L_{p,n}$}
\psfrag{K}{\footnotesize $K$}
\includegraphics[width=12cm]{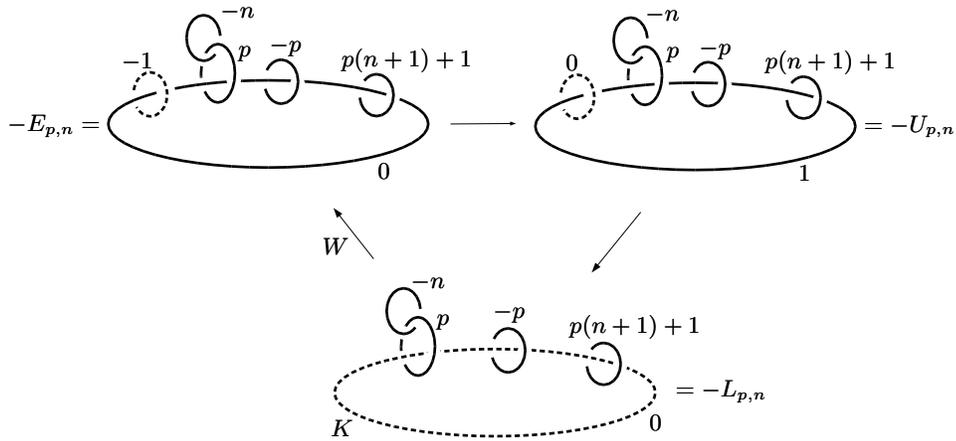}
%\put(0,0){\tiny\grid(130,110)(5,5)[0,0]}
\end{center}
\caption{Manifolds and cobordisms in a related  surgery triangle}
\label{f:reltri}
\end{figure}

Let us consider the triangle of 3--manifolds and cobordisms
given by Figure~\ref{f:reltri}.

\begin{prop}
The 3--manifold $-U_{p,n}$ is an $L$--space.
\end{prop}

\begin{proof}
Kirby calculus, as in the proof of Proposition~\ref{p:slides}, shows
that $U_{p,n}$ is diffeomorphic to $S^3_r(T_{p,pn+1})$, with
\[
r=p^2n+p+1+\frac{1}{p(n+1)}. 
\]
Since the above $r$ is greater than $2g_s(T_{p,pn+1})-1=p^2n-pn-1$,
by~\cite[Proposition~4.1]{LSuj} $U_{p,n}$ is an $L$--space.
\end{proof}

The exact triangle on Ozsv\'ath--Szab\'o homologies induced by the
surgery triangle of Figure~\ref{f:reltri} has the following shape:
\[
\begin{graph}(6,2)
\graphlinecolour{1}\grapharrowtype{2}
\textnode {A}(1,1.5){$\hf (-E_{p,n})$}
\textnode {B}(5, 1.5){$\hf (-U_{p,n})$}
\textnode {C}(3, 0){$\hf (-L_{p,n})$}
\diredge {A}{B}[\graphlinecolour{0}]
\diredge {B}{C}[\graphlinecolour{0}]
\diredge {C}{A}[\graphlinecolour{0}]
\freetext (3,1.8){$F'$}
\freetext (4.6,0.6){$H'$}
\freetext (1.4,0.6){$G_W$}
\end{graph}
\]
Simple computation shows that
\[
h_U:=\vert H_1(U_{p,n}; \Z )\vert =p^3n(n+1)+p(p+1)(n+1)+1. 
\]

Since $h_U$ is odd, the 3--manifold $-U_{p,n}$ supports a
unique spin structure, which will be denoted by $\t_U$. In analogy to
Equation~\eqref{e:dsum}, there is a direct sum decompositions
\begin{eqnarray}\label{e:decompu}
\hf(-U_{p,n}) = \langle \t_U \rangle\oplus S\oplus\overline{S}.
\end{eqnarray}

\begin{cor}
The map $F'$ in the above triangle is 0. Therefore
$H'$ is injective and $G_W$ is surjective.
\end{cor}
\begin{proof}
  Since all the manifolds involved are $L$--spaces, the argument boils down to
  the simple observation that $h_L=h_E+h_U$, cf. also the proof of
  Proposition~\ref{p:van}.
\end{proof}

\begin{lem}\label{l:centspinc}
The $\t_E$--component of the element $G_W(\t_W)\in \hf (-E_{p,n})$
is nonzero.
\end{lem}

\begin{proof}
  Notice first that, since $\t_V$ does not extend to $W$ as a spin
  structure, by Lemma~\ref{l:ext} the $\t_E$--component of $G_W(\t_V)$
  is zero.  Arguing by contradiction, suppose now that the
  $\t_E$--component of $G_W(\t_W)$ is also zero. Suppose that $G_W
  (\t_V)=x_V + {\overline {x_V}}$ and $G_W (\t_W)=x_W + {\overline
  {x_W}}$ with $x_V, x_W\in T$.  

  Since $G_W$ is onto, there exist
  elements $l_V, l_W\in\hf(-L_{p,n})$ such that 
  \[
  G_W(l_V)=x_V\quad\text{and}\quad G_W(l_W)=x_W.  
  \] 
  Therefore,  
  \[
  G_W(\t_V + l_V +\overline{l_V}) = 0\quad\text{and}\quad G_W(\t_W+l_W
  +\overline{l_W})=0.  
  \] 
  By exactness, this implies the existence of 
  $u_V, u_W\in\hf(-U_{p,n})$
  such that 
  \[ H'(u_V)=\t_V+ l_V +\overline{l_V}\quad\text{and}\quad
  H'(u_W)=\t_W+l_W +\overline{l_W}.  
  \] 
  Since $H'$ is injective, we
  have that $u_V$ and $u_W$ are both fixed under conjugation. Then,
  one of $u_V$, $u_W$ or $u_V+u_W$ belongs to $S\oplus\overline{S}$
  and is therefore of the form $s+\overline{s}$ for some $s\in S$. But
  for any $s\in S$ we have $H'(s+\overline{s})\in
  C\oplus\overline{C}$, so one of $t_V+l_V+{\overline {l_V}}$,
  $t_W+l_W+{\overline {l_W}}$ 
  or their sum belongs to $C\oplus {\overline {C}}$,
  which is clearly impossible. This
  contradiction proves the lemma.
\end{proof}
  
The following is the most important result of this subsection

\begin{prop}\label{p:comp}
We have
\[
\deg(\t_E)=\deg(\t_W)+\frac{1}{4}.
\]
\end{prop}

\begin{proof}
  By Lemma~\ref{l:centspinc} the element $G_W(\t_W)$ has nontrivial
  $\t_E$--coordinate, therefore there are spin$^c$ structures $\s_i$ on $W$
  such that $G_{W,\s_i}(\t_W)=\t_E$. By the conjugation invariance we have
  that $G_{W,\s_i}(\t_W)= G_{W, {\overline {\s_i}}}(\t_W)$. Since we use mod 2
  coefficients, this shows that there are an odd number of $\s_i$'s with the
  above property, and therefore there exists a spin structure $\s$ on $W$ with
  the property that $G_{W,\s}(\t_W)=\t_E$. An argument similar to the one
  given in Lemma~\ref{l:signatures} shows that $W$ is negative definite. Since
  for a spin structure $c_1(\s)=0$, the degree shift formula implies the
  result.
\end{proof}

\begin{proof}[Proof of Theorem~\ref{t:T34}]
Recall that there is an element $a\in \hf (-S_{p,n})$ satisfying the
equation $G_V(\t_W)=a +{\overline {a}}$. Express $a$ as a sum of
homogeneous elements. Since by Lemma~\ref{l:tE-comp} the
$\t_E$--component of $F(a)$ is nonzero, $a$ has a homogeneous
component $a_1$ with the same property. By the degree--shift formula,
Lemmas~\ref{l:c1V} and~\ref{l:c1X} immediately imply (with $\vert k \vert =
\vert l \vert =1$) that
\begin{equation}\label{e:first}
gr(\t_E)-\frac{1}{4} (-\frac{h_E}{h_S}+1)\leq
gr(a_1)\leq gr (\t_W)+\frac{1}{4} (-\frac{h_L}{h_S}+1).
\end{equation}
But since $h_S=h_L+h_E$, by Proposition~\ref{p:comp} the inequalities
of Equation~\eqref{e:first} must in fact be equalities. This shows
that the spin$^c$ structure corresponding to $a_1$ is the restriction
of a spin$^c$ structure $\s$ as in Lemma~\ref{l:c1V} with $k=\pm
1$. Consequently, $a_1\in \hf (-S_{p,n}, \t )$ with $c_1 (\t ) = \pm
PD (\mu _d)$ in the basis of homologies given by
Figure~\ref{f:traf}. According to Lemma~\ref{l:spinccomp}, either
$a_1$ or ${\overline {a_1}}$ belongs to the same summand
$\hf(-S_{p,n},\t)$ as $c(S_{p,n},\xi_{p,n})$. Therefore, since
$-S_{p,n}$ is an $L$--space, $c(S_{p,n},\xi_{p,n})$ is equal to either
$a_1$ or $\overline{a_1}$.  But
$F(\overline{a_1})=\overline{F(a_1)}$. Therefore,
\[
c(E_{p,n}, \zeta_{p,n}) = F_{-X}(c(S_{p,n}, \xi_{p,n})) 
\]
has nonzero $\t_E$--component, and therefore it coincides with $\t_E$. This
fact implies that $\zeta _{p,n}$ is a tight, positive contact structure on
$E_{p,n}$, concluding the proof.
\end{proof}

\end{document}